\documentclass{amsart}
\input{preamble.tex}

\graphicspath{{Fig/}}

\begin{document}

% \DOI{DOI HERE}
% \copyrightyear{2023}
% \vol{00}
% \pubyear{2023}
% \access{Advance Access Publication Date: Day Month Year}
% \appnotes{Paper}
% \copyrightstatement{Published by Oxford University Press on behalf of the Institute of Mathematics and its Applications. All rights reserved.}
% \firstpage{1}

% \setcounter{section}{1}

\title[Geometric integration for smooth optimisation]{A geometric integration approach to smooth optimisation: Foundations of the discrete gradient method}

\author[M. J. Ehrhardt \emph{et al.}]{Matthias J. Ehrhardt}
\address[MJE]{Department of Mathematical Sciences, University of Bath, North Rd, Claverton Down, Bath, BA2 7AY, UK}

\author[]{Erlend S. Riis}
\address[ESR]{Department of Applied Mathematics and Theoretical Physics, University of Cambridge, Wilberforce Rd, Cambridge, CB3 0WA, UK}
\email{erlendsriis@gmail.com}

\author[]{Torbj{\o}rn Ringholm}
\address[TR]{Department of Mathematical Sciences, Norwegian University of Science and Technology, Alfred Getz' vei 1, 7034 Trondheim, Norway}

\author[]{Carola-Bibiane Sch{\"o}nlieb}
\address[CBS]{Department of Applied Mathematics and Theoretical Physics, University of Cambridge, Wilberforce Rd, Cambridge, CB3 0WA, UK}

\begin{abstract}
	Discrete gradient methods are geometric integration techniques that can preserve the dissipative structure of gradient flows. Due to the monotonic decay of the function values, they are well suited for general convex and nonconvex optimisation problems. Both zero- and first-order algorithms can be derived from the discrete gradient method by selecting different discrete gradients. In this paper, we present a thorough analysis of the discrete gradient method for optimisation which provides a solid theoretical foundation. We show that the discrete gradient method is well-posed by proving the existence of iterates for any positive time step, as well as uniqueness in some cases, and propose an efficient method for solving the associated discrete gradient equation. Moreover, we establish an \(\rm{O}(1/k)\) convergence rate for convex objectives and prove linear convergence if instead the Polyak--Łojasiewicz inequality is satisfied. The analysis is carried out for three discrete gradients—the Gonzalez discrete gradient, the mean value discrete gradient, and the Itoh--Abe discrete gradient—as well as for a randomised Itoh--Abe method. Our theoretical results are illustrated with a variety of numerical experiments, and we furthermore demonstrate that the methods are robust with respect to stiffness.
\end{abstract}

\maketitle

\section{Introduction}\label{intro}

Discrete gradients are tools from geometric integration for numerically solving first-order systems of ordinary differential equations (ODEs), while ensuring that certain structures of the continuous system---specifically energy conservation and dissipation, and Lyapunov functions---are preserved in the numerical solution. Energy dissipation refers to the monotonic decrease in value of the objective function over time.

The use of discrete gradient methods to solve optimisation problems has gained increasing attention in recent years \citep{cel18,gri17,rii21,rin18}, due to their preservation of dissipative structures of ODEs such as gradient flows. This means that the associated iterative scheme monotonically decreases the objective function for all positive time steps, at a rate analogous to that of gradient flow.

In this paper, we consider the unconstrained optimisation problem
\begin{equation}\label{eq:optimisation_problem}
	\min_{x \in \RR^n} V(x),
\end{equation}
where the function \(V: \RR^n \to \RR\) is continuously differentiable. For an initial guess \(x^0 \in \RR^n\) and \(k \in \NN\), the discrete gradient method is of the form
\begin{equation}\label{eq:dg_method}
	x^{k+1} = x^k - \tau_k \overline{\nabla} V(x^k, x^{k+1}),
\end{equation}
where \(\tau_k > 0\) is the time step, and \(\overline{\nabla} V\) is the discrete gradient, defined as follows.
\begin{definition}[Discrete gradient]\label{defn:dg}
	Let \(V\) be a continuously differentiable function. A \emph{discrete gradient} is a continuous map \(\overline{\nabla} V: \RR^n \times \RR^n \to \RR^n\) such that for all \(x, y \in \RR^n\),
	\begin{alignat}{2}\label{eq:mean_value}
		\inner{\overline{\nabla} V(x,y), y-x}           & = V(y) - V(x) & \qquad & \mbox{(Mean value property)},  \\
		\label{eq:consistency} \overline{\nabla} V(x,x) & = \nabla V(x) & \qquad & \mbox{(Consistency property)}.
	\end{alignat}
\end{definition}
The background for discrete gradients and gradient flows for optimisation is given in \secref{sec:dg}.

There are several aspects of discrete gradient methods which make them attractive for optimisation. Their structure-preserving properties lead to schemes that are unconditionally stable, i.e.\ converging to a stationary point for arbitrary time steps. This may be particularly beneficial for stiff problems, where explicit schemes require prohibitively small time steps, due to local, rapid variations in the objective function gradient. We demonstrate this numerically in \secref{sec:TV}. Furthermore, the Itoh--Abe discrete gradient~\eqref{eq:itoh_abe_dg} is derivative-free, enabling derivative-free optimisation\footnote{Not to be confused with the derivative-free discrete gradient method proposed by~\cite{bag08} which uses a different notion of a discrete gradient.} e.g.\ for nonsmooth, nonconvex functions and \emph{black-box} problems \citep{rii21}, and for functions whose gradients are expensive to compute \citep{rin18}. Beyond Euclidean gradient flow, discrete gradients can also preserve the dissipative structure for inverse scale space and Bregman distances \citep{ben20}.

\subsection{Contributions and structure}

The aim of this paper is to give a comprehensive optimisation analysis of discrete gradient methods. While having a solid theoretical foundation in geometric integration, the methods have only recently been explored as optimisation schemes, leaving several gaps in our understanding of their properties.

To this end, we address several theoretical questions. Namely, we prove well-posedness of the implicit update equation~\eqref{eq:dg_method}, propose efficient and stable methods for solving the update equation, and obtain convergence rates of the methods for different classes of objective functions. Furthermore, we provide various numerical examples to see how the methods perform in practice.

In \secref{sec:dg} we define discrete gradients and introduce the four discrete gradient methods considered in this paper. In \secref{sec:existence}, we prove that the discrete gradient equation (the update formula)~\eqref{eq:dg_method} is well-posed, meaning that for any time step \(\tau_k > 0\) and \(x^k \in \RR^n\), a solution \(x^{k+1}\) exists, under mild assumptions on \(V\). Using Brouwer's fixed point theorem, this is the first existence result for the discrete gradient equation without a bound on the time step. In \secref{sec:fixed_point}, we propose an efficient and stable method for solving the discrete gradient equation and prove convergence guarantees, building on the ideas of~\cite{nor14}.

In \secref{sec:analysis_time_step}, we analyse the dependence of the iterates on the choice of time step, and obtain estimates for preferable time steps in the cases of \(L\)-smoothness and strong convexity. In \secref{sec:convergence_rate}, we establish convergence rates for convex functions with Lipschitz continuous gradients, and for functions that satisfy the Polyak--{\L}ojasiewicz (P{\L}) inequality \citep{kar16}. In \secref{sec:preconditioned}, we briefly discuss preconditioned discrete gradient methods.

In \secref{sec:numerical}, we present numerical results for several test problems, and a numerical comparison of different numerical solvers for the discrete gradient equation~\eqref{eq:dg_method}. We conclude and present an outlook for future work in \secref{sec:conclusion}.

We emphasise that the majority of these results hold for nonconvex functions. Our contributions to the foundations of the discrete gradient method opens the door for future applications and research on discrete gradient methods for optimisation, with a deeper understanding of their numerical properties.

\subsection{Related work}

We list some applications of discrete gradient methods for optimisation. \cite{gri17} propose applying them to solve variational problems in image analysis and prove convergence to stationary points for continously differentiable functions. \cite{rin18} apply the Itoh--Abe discrete gradient method to nonconvex imaging problems regularised with Euler's elastica. \cite{miy18} point out the equivalence between the Itoh--Abe discrete gradient method for quadratic problems and the well-known Gauss-Seidel and successive-over-relaxation (SOR) methods. % chktex 2 chktex 12

Several recent works look at discrete gradient methods in other optimisation settings. Concerning nonsmooth, derivative-free optimisation, \cite{rii21} study the Itoh--Abe discrete gradient method for solving nonconvex, nonsmooth functions, proving converge to a set of stationary points in the Clarke subdifferential framework. \cite{cel18} extend the Itoh--Abe discrete gradient method for optimisation on Riemannian manifolds. \cite{her15} combine a discrete gradient method with Hopfield networks to preserve a Lyapunov function for optimisation problems. \cite{ben20} apply the discrete gradient method to inverse scale space flow for sparse optimisation. % chktex 2 chktex 12

More generally, there is a wide range of research that studies connections between optimisation schemes and systems of ODEs. \cite{su16} and \cite{wib16} study second-order ODEs that can be seen as continuous-time limits of Nesterov acceleration schemes. In the former case, they shed light on the dynamics of these schemes, e.g.\ the oscillatory behaviour, by interpreting the ODEs as damping systems, while in the latter case, they present a family of \emph{Bregman Lagrangian functionals} which generate the original and new acceleration schemes. Furthermore, they demonstrate that the choice of ODE discretisation method is central for whether the acceleration phenomena is retained in the iterative scheme. \cite{att00} consider the second-order heavy ball with friction dynamical system for convex optimisation problems in Hilbert spaces. Several works have contributed to the setting of numerical analysis of acceleration methods which bridges continous- and discrete-time dynamics. \cite{wil21} approach this from the perspective of Lyapunov theory, presenting Lyapunov functions that account for continuous- and discrete-time. \cite{bet18} present a framework of \emph{sympletic optimisation}, considering perspectives of Hamiltonian dynamics and symplectic structure-preserving methods. In a similar vein, recent papers by \cite{mad18} and \cite{fra20} study conformal Hamiltonian systems, with the former using information about the objective function's convex conjugate to obtain stronger convergence rates, and the latter considering structure-preserving numerical methods and their relation to some iterative schemes. % chktex 2 chktex 12

\subsection{Notation and preliminaries}

We denote by \(S^{n-1}\) the unit sphere \(\set{x \in \RR^n \, : \, \|x\| = 1}\). The line segment between two points is defined as \([x,y] := \set{\lambda x + (1-\lambda)y \; : \; \lambda \in [0,1]}\). The diameter of a set \(U \subset \RR^n\) is defined as \(\diam(U) := \sup_{x, y \in U} \|x-y\|\).

In this paper, we consider both deterministic schemes and stochastic schemes. For the stochastic schemes, there is a random distribution \(\Xi\) on \(S^{n-1}\) such that each iterate \(x^k\) depends on a descent direction \(d^k\) which is independently
drawn from \(\Xi\). We denote by \(\xi^k\) the joint distribution of \({(d^i)}^k_{i=1}\), and by \(V_{k}\) the expectation of \(V(x^{k})\) with respect to \(\xi^{k}\),
\begin{equation}\label{eq:expectation}
	V_{k} := \EE[ V(x^{k}) ].
\end{equation}
We also use the shorthand notation \(\EE_d[V(x)]\) for \(\EE_{d \sim \Xi}[V(x)]\). To unify notation for all the methods in this paper, we will write \(V_{k}\) instead of \(V(x^{k})\) for the deterministic methods as well.

Throughout the paper, we will consider two classes of functions, \emph{\(L\)-smooth} and \emph{\(\mu\)-convex} functions. We here provide definitions and some basic properties.
\begin{definition}[\(L\)-smooth]\label{defn:smooth}
	A function \(V: \RR^n \to \RR\) is \emph{\(L\)-smooth} for \(L > 0\) if its gradient is Lipschitz continuous with Lipschitz constant \(L\), i.e.\ if for all \(x,y \in \RR^n\),
	\[
		\|\nabla V(x) - \nabla V(y)\| \leq L\|x-y\|.
	\]
\end{definition}
We state some properties of \(L\)-smooth functions.
\begin{proposition}\label{prop:smooth}
	If \(V:\RR^n \to \RR\) is \(L\)-smooth, then for all \(x, y \in \RR^n\), the following holds.
	\begin{enumerate}[(i) ]
		\item \(V(y) - V(x) \leq \inner{\nabla V(x), y-x} + \dfrac{L}{2}\|y-x\|^2\).
		\item \(V(\lambda x + (1-\lambda) y) \geq  \lambda V(x) + (1-\lambda) V(y) - \dfrac{\lambda (1-\lambda) L}{2} \|x-y\|^2\) for all \(\lambda \in [0,1]\).
	\end{enumerate}
\end{proposition}

\begin{proof}
	Property \emph{(i)}. \citep[Proposition A.24]{ber03}. % chktex 12

	Property \emph{(ii)}. It follows from property \emph{(i)} that the function \(x \mapsto \frac{L}{2}\|x\|^2 - V(x)\) is convex, which in turn yields the desired inequality.
\end{proof}

\begin{definition}[$\mu$-convex]\label{defn:strongly_convex}  % chktex 46
	A function \(V: \RR^n \to \RR\) is \emph{\(\mu\)-convex} for \(\mu \geq 0\) if either of the following (equivalent) conditions hold.
	\begin{enumerate}[(i) ]
		\item The function \(V(\cdot) - \dfrac{\mu}{2} \|\cdot\|^2\) is convex.
		\item \(V\del{\lambda x + (1-\lambda) y} \leq \lambda V(x) + (1-\lambda) V(y) - \dfrac{\mu}{2} \lambda (1-\lambda) \|x-y\|^2\) for all \(x,y\) in \(\RR^n\) and \(\lambda \in [0,1]\).
	\end{enumerate}
	When \(\mu > 0\), we say that \(V\) is \emph{strongly convex}.
\end{definition}

\section{Discrete gradient methods}\label{sec:dg}

\subsection{The discrete gradient method and gradient flow}

We motivate the use of discrete gradients by considering the gradient flow of \(V\),
\begin{equation}\label{eq:gradient_flow}
	\dot{x} = -\nabla V(x), \qquad x(0) = x_0 \in \RR^n,
\end{equation}
where \(\dot{x}\) denotes the derivative of \(x\) with respect to time. This system is fundamental to optimisation, and underpins many gradient-based schemes. Applying the chain rule, we obtain
\begin{equation}\label{eq:gradient_flow_dissipation}
	\dod{}{t}V(x(t)) = \inner{\nabla V(x(t)), \dot{x}(t)} = - \|\nabla V(x(t))\|^2 = - \|\dot{x}(t)\|^2 \leq 0.
\end{equation}
The gradient flow has an energy dissipative structure, since the function value \(V(x(t))\) decreases monotonically along any solution \(x(t)\) to~\eqref{eq:gradient_flow}. Furthermore, the rate of dissipation is given in terms of the norm of \(\nabla V\) or equivalently the norm of \(\dot{x}\).

In geometric integration, one studies methods for numerically solving ODEs while preserving certain structures of the continuous system---see~\cite{hai06, mcl01} for an introduction. Discrete gradients are tools for solving first-order ODEs that preserve energy conservation laws, dissipation laws, and Lyapunov functions \citep{gon96, ito88, mcl99, qui96}.

For a sequence of strictly positive time steps \({(\tau_{k})}_{k \in \NN}\) and a starting point \(x^0 \in \RR^n\), the discrete gradient method applied to~\eqref{eq:gradient_flow} is given by~\eqref{eq:dg_method}.
This scheme preserves the dissipative structure of gradient flow, as we see by applying~\eqref{eq:mean_value},
\begin{equation}\label{eq:dissipation}
	\begin{aligned}
		V(x^{k+1}) - V(x^{k}) & =  \inner{\overline{\nabla}V(x^k, x^{k+1}), x^{k+1} - x^{k}} \\ & = - \tau_k \|\overline{\nabla}V(x^k, x^{k+1})\|^2 = - \frac{1}{\tau_{k}} \|x^{k+1} - x^k\|^2.
	\end{aligned}
\end{equation}
Similarly to the dissipation law~\eqref{eq:gradient_flow_dissipation} of gradient flow, the decrease of the objective function value is given in terms of the norm of both the step \(x^{k+1} - x^k\) and of the discrete gradient.

Throughout the paper, we assume there are bounds \(\tau_{\max} \geq \tau_{\min} > 0\) such that
\begin{equation}\label{eq:time_step}
	{(\tau_k)}_{k\in\NN} \subset [\tau_{\min}, \tau_{\max}].
\end{equation}
No restrictions are required on these bounds. \cite{gri17} prove that if \(V\) is continuously differentiable and coercive---the latter meaning that the level set \(\set {x \in \RR^n \; : \; V(x) \leq M}\) is bounded for each \(M \in \RR\)---and if \({(\tau_k)}_{k \in \NN}\) satisfy~\eqref{eq:time_step}, then the iterates \({(x^{k})}_{k \in \NN}\) of~\eqref{eq:dg_method} converge to a set of stationary points, i.e.\ points \(x^* \in \RR^n\) such that  \(\nabla V(x^*) = 0\). % chktex 2  % chktex 12

We may compare the discrete gradient method to explicit gradient descent, \(x^{k+1} = x^k - \tau_k \nabla V(x^{k})\). Unlike discrete gradient methods, gradient descent only decreases the objective function value for sufficiently small time steps \(\tau_k\). To ensure decrease and convergence for this scheme, the time steps must be restricted based on estimates of the smoothness of the gradient of \(V\), which might be unavailable, or lead to prohibitively small time steps. Conversely, implicit gradient descent, or the proximal point method~\citep{cha16}, is unconditionally stable with respect to the time step, but does not necessarily exhibit structure preservation~\citep{gri17}.

One may also consider other numerical integration methods, such as implicit Runge-Kutta methods, where energy dissipation is ensured under mild time step restrictions~\citep{hai13}, and explicit stabilised methods for solving strongly convex problems~\citep{eft21}.

\subsection{Four discrete gradient methods}

We now introduce the four discrete gradients considered in this paper.
% TO DO: Consider switching to enumerate environment
1.\ \emph{The Gonzalez discrete gradient} \citep{gon96}  is given by
\[
	\overline{\nabla} V(x,y) = \nabla V\del{\frac{x+y}{2}} + \dfrac{V(y) - V(x) - \inner{\nabla V(\frac{x+y}{2}), y - x}}{\|x-y\|^2} (y-x), \quad x \neq y.
\]

2.\ \emph{The mean value discrete gradient} \citep{har83}, used for example in the average vector field method~\citep{cel12}, is given by
\[
	\overline{\nabla} V(x,y) = \int_0^1 \! \nabla V\del{(1-s)x + sy} \dif s.
\]

3.\ \emph{The Itoh--Abe discrete gradient} \citep{ito88} (also known as the coordinate increment discrete gradient) is given by
\begin{equation}\label{eq:itoh_abe_dg}
	\overline{\nabla} V(x,y) = \begin{pmatrix}
		\frac{V\del{y_1, x_2, \ldots, x_n} - V(x)}{y_1 - x_1}                              \\
		\frac{V\del{y_1, y_2, x_3, \ldots, x_n} - V\del{y_1, x_2, \ldots, x_n}}{y_2 - x_2} \\
		\vdots                                                                             \\
		\frac{V(y) - V\del{y_1, \ldots, y_{n-1}, x_n}}{y_n - x_n}
	\end{pmatrix},
\end{equation}
where \(0/0\) is interpreted as \(\partial_i V(x)\).

While the first two discrete gradients are gradient-based and can be seen as approximations to the midpoint gradient \(\nabla V(\tfrac{x+y}{2})\), the Itoh--Abe discrete gradient is derivative-free, and is evaluated by computing successive coordinate-wise difference quotients. As a result, the corresponding implicit equation~\eqref{eq:dg_method} uncouples, and can be treated as solving a series of scalar equations
\begin{equation*}
	x^{k+1}_i = x^k_i - \tau_k \frac{V(x^{k+1}_1, \ldots, x^{k+1}_i, x^k_{i+1}, \ldots, x_n^{k})- V(x^{k+1}_1, \ldots, x^{k+1}_{i-1},x^k_i,\ldots, x_n^{k})}{x^{k+1}_i - x^k_{i}},
\end{equation*}
for \(k=1,\ldots,n\). In an optimisation setting, the Itoh--Abe discrete gradient is often preferable to the others, as it is relatively computationally inexpensive.

4.\ \emph{The Randomised Itoh--Abe method} \citep{rii21} is an extension of the Itoh--Abe discrete gradient method, where the directions of descent are randomly chosen. Given a sequence of directions \({(d^{k})}_{k \in \NN} \subset S^{n-1}\) drawn independently from a random distribution \(\Xi\), we solve
\begin{equation*}\label{eq:random_IA}
	x^{k+1} \mapsto x^k - \tau_k \alpha_k d^{k+1}, \quad  \mbox{where } \alpha_k \neq 0 \mbox{ solves } \quad \alpha_k = - \frac{V(x^k - \tau_k \alpha_k d^{k+1}) - V(x^{k})}{\tau_k \alpha_{k}},
\end{equation*}
where \(x^{k+1} = x^k\) is considered a solution whenever \(\inner{\nabla V(x^{k}), d^{k+1}} = 0\).

This scheme generalises the Itoh--Abe discrete gradient method, in that the methods are equivalent if \({(d^{k})}_{k \in \NN}\) cycle through the standard coordinates with the rule \(d^k = e^{[(k-1) \Mod n] + 1}\),  \(k = 1, 2, \ldots\).	However, the computational effort of one iterate of the Itoh--Abe discrete gradient method equals \(n\) steps of the randomised method, so their efficiency should be judged accordingly.

While the randomised Itoh--Abe discrete gradient method does not retain the discrete gradient structure of the Itoh--Abe discrete gradient, it retains a dissipative structure akin to~\eqref{eq:gradient_flow_dissipation}.
\begin{align}\label{eq:dissipation_IA}
	V(x^{k+1}) - V(x^{k}) & = -\tau_k\del{\frac{V(x^{k+1}) - V(x^{k})}{\|x^{k+1} - x^k\|}}^2 = -\frac{1}{\tau_{k}} \|x^{k+1} - x^k\|^2.
\end{align}

We also define the constant
\begin{equation}\label{eq:angle_bound}
	\zeta := \min_{e \in S^{n-1}} \EE_d [ \inner{d, e}^2 ],
\end{equation}
and assume that \(\Xi\) is such that \(\zeta > 0\). For example, for the uniform random distribution on both \(S^{n-1}\) and on the standard coordinates \({(e^{i})}_{i=1}^n\), we have \(\zeta = 1/n\). See \citep[Table 4.1]{sti14} for estimates of~\eqref{eq:angle_bound} for these cases and others.

The motivation for introducing this randomised extension of the Itoh--Abe method is, first, to tie in discrete gradient methods with other optimisation methods such as stochastic coordinate descent \citep{fer15, qu15, wri15} and random pursuit \citep{nes17, sti14}, and, second, because this method extends to the nonsmooth, nonconvex setting \citep{rii21}.

\section{Existence of solutions to the discrete gradient steps}\label{sec:existence}

In this section, we prove that the discrete gradient equation
\begin{equation}\label{eq:dg}
	y = x - \tau \overline{\nabla} V(x, y).
\end{equation}
admits a solution \(y\), for all \(\tau > 0\) and \(x \in \RR^n\), under mild assumptions on \(V\) and \(\overline{\nabla} V\).

To the authors' knowledge, the following result is the first without a restriction on time steps. \cite{nor14} provide an existence and uniqueness result for small time steps for a large class of discrete gradients, via the Banach fixed point theorem. Furthermore, the existence of a solution for the Gonzalez discrete gradient is established for sufficiently small time steps via the implicit function theorem in \citep[Theorem 8.5.4]{stu96}. % chktex 2 chktex 12

We use the following set notation. For \(\delta > 0\), the closed ball of radius \(\delta\) about \(x\) is defined as \(\overline{B}_\delta(x) := \set{y \in \RR^n \; : \; \|y - x\| \leq \delta}\). For a set \(U \subset \RR^n\), its \(\delta\)-thickening is the set defined as \({U_\delta := \set{x \in \RR^n \; : \; \dist(U,x) \leq \delta}}\). The convex hull of \(U\) is given by \(\co(U)\).

We make two assumptions for the discrete gradient, namely that boundedness of the gradient implies boundedness of the discrete gradient, and that if two functions coincide on an open set, their discrete gradients also coincide.
\begin{assumption}\label{ass:dg}
	There is a constant \(C_n\) that depends on the discrete gradient but is independent of \(V\), and a continuous, nondecreasing function \(\delta: [0, \infty] \to [0, \infty]\), where \(\delta(0) = 0\) and \(\delta(\infty) := \lim_{r \to \infty} \delta(r)\), such that the following holds.

	For any \(V \in C^1(\RR^n ; \RR)\) and any convex set \(U \subset \RR^n\) with nonempty interior, the two following properties are satisfied.
	\begin{enumerate}[(i) ]
		\item If \(\|\nabla V(x)\| \leq K\) for all \(x \in U_{\delta\del{\diam(U)}}\), then \(\|\overline{\nabla} V(x,y)\| \leq C_n K\) for all \(x, y \in U\).
		\item If \(W\) is another continuously differentiable function such that \(V(x) = W(x)\) for all \(x \in U_{\delta\del{\diam(U)}}\), then \(\overline{\nabla} V(x,y) = \overline{\nabla} W(x,y)\) for all \(x,y \in U\).
	\end{enumerate}
\end{assumption}
\begin{remark}
	It is straightforward to verify that if a finite collection of discrete gradient constructions satisfy these assumptions, then their convex combinations satisfy them too.
\end{remark}

The following result, which is proved in \appref{sec:dg_bounds}, shows that the discrete gradients considered in this paper satisfy the above assumption.

\begin{lemma}\label{lem:dg_assumption}
	The three discrete gradients satisfy \assref{ass:dg} as follows.
	\begin{enumerate}
		\item For the Gonzalez discrete gradient, \(C_n = \sqrt{2}\) and \(\delta \equiv 0\).
		\item For the mean value discrete gradient, \(C_n = 1\) and \(\delta \equiv 0\).
		\item For the Itoh--Abe discrete gradient, \(C_n = \sqrt{n}\) and \(\delta(r) = r\).
	\end{enumerate}
\end{lemma}

\begin{remark}
	In practice, the Gonzalez, Itoh--Abe, and mean value discrete gradients account for the vast majority of applications and theoretical studies \citep{dah11,gri17,mcl99,qui96}. Thus, while \assref{ass:dg} may not hold for all constructions of discrete gradients, we consider them to be adequate for most purposes.
\end{remark}

The existence proof is based on the  Brouwer fixed point theorem \citep{bro11}.
\begin{proposition}[Brouwer fixed point theorem]\label{prop:brouwer}
	Let \(U \subset \RR^n\) be a convex, compact set and \(g: U \to U\) a continuous function. Then \(g\) has a fixed point in \(U\).
\end{proposition}

We proceed to state the existence theorem.
\begin{theorem}[Discrete gradient existence theorem]\label{thm:existence}
	Suppose \(V\) is continously differentiable and that \(\overline{\nabla}\) satisfies \assref{ass:dg}. Then there exists a solution \(y\) to~\eqref{eq:dg} for any \(\tau>0\) and \(x \in \RR^n\), if \(V\) satisfies \textbf{any} of the following properties.
	\begin{enumerate}[(i) ]
		\item The gradient of \(V\) is uniformly bounded.
		\item \(V\) is coercive.
		\item Both \(V\) and the gradient of \(V\) are uniformly bounded on \(\co(\{y \; : \; V(y) \leq V(x)\})\) (the bounds may depend on \(x\)), and \(\delta \equiv 0\) in \assref{ass:dg}.
	\end{enumerate}
\end{theorem}

\begin{proof}
	Part \emph{(i)}. We want to show that the function \(g(y) = x - \tau \overline{\nabla}V(x,y)\) has a fixed point, \(y = g(y)\). There is \(K > 0\) such that \(\|\nabla V(y) \| \leq K\) for all \(y \in \RR^n\). Therefore, by \assref{ass:dg}, \(\|\overline{\nabla} V(x,y)\| \leq C_n K\) for all \(y \in \RR^n\). This implies that \(g(y) \in \overline{B}_{\tau C_n K}(x)\) for all \(y \in \RR^n\). Specifically, \(g\) maps \(\overline{B}_{\tau C_n K}(x)\) into itself. As \(g\) is continuous, it follows from \propref{prop:brouwer} that there exists a point \(y \in \overline{B}_{\tau C_n K}(x)\) such that \(g(y) = y\), and we are done.

	Part \emph{(ii)}. Let \(\sigma >  0\), \(U=\co(\set{y \; : \; V(y) \leq V(x)})\), and write \(\delta = \delta(\diam(U))\). Since \(V\) is coercive, \(U_\delta\) and \(U_{\delta+\sigma}\) are bounded. By standard arguments \citep[Corollary 2.5]{nes03}, there exists a cutoff function \(\varphi \in C_c^\infty(\RR^n; [0,1])\) such that \(\varphi\mid_{U_\delta} \equiv 1 \) and \(\varphi \mid_{\RR^n \setminus U_{\delta+\sigma}} \equiv 0\).
	We then define the function \[W(y) := \varphi(y)\del{V(y) -  V(x)} + V(x),\] which is continuously differentiable and \(\supp(\nabla W) \subset U_{\delta + \sigma}\). Therefore, \(\nabla W\) is uniformly bounded, so by part \emph{(i)} there is \(y\) such that
	\[
		y = x - \tau \overline{\nabla} W(x,y).
	\]
	By~\eqref{eq:dissipation}, \(W(y) < W(x)\), which, by construction of \(W\), implies that \(V(y) < V(x)\) and hence that \(y \in U\). Lastly, since \(V\) and \(W\) coincide on \(U_\delta\), and \(x\) and \(y\) both belong to \(U\), it follows from \assref{ass:dg} \emph{(ii)} that \(\overline{\nabla} V(x,y) = \overline{\nabla} W(x,y)\). Hence a solution \(y = x - \tau \overline{\nabla} V(x,y)\) exists.

	Part \emph{(iii)}. Set \(U = \co\set{y \; : \; V(y) \leq V(x)}\) and \(M = \sup_{y \in U} V(y)\). Furthermore let \(\eta > 0\) and set \(K = \sup\set{ \|\nabla V(y)\| \; : \; V(y) \leq M + \eta}\) and \(F = \set{y \; : \; V(y) \geq M + \eta}\). For any \(z \in F\) and \(y \in U\), there is \(\overline{z} \in [z,y]\) such that \(V(\overline{z}) = M+\eta\). The mean value theorem (MVT) \citep[Equation A.55]{noc99} and the boundedness of \(\nabla V\) imply that there is \(\lambda \in (0,1)\) such that
	\[
		\eta \leq |V(y) - V(\overline{z})| = |\inner{\nabla V(\lambda y + (1-\lambda)\overline{z}), y - \overline{z}}| \leq K \|y - \overline{z}\| \leq K \|y - z\|.
	\]
	Therefore, for all \(y \in U\) and \(z \in F\), one has \(\|y - z\| \geq \eta/K\). As in the previous case, there exists a cutoff function \(\varphi \in C^\infty(\RR^n; [0,1])\) with uniformly bounded gradient, such that \(\varphi \mid_{U} \equiv 1\) and \(\varphi\mid_{F} \equiv 0\).
	To verify this, we may consider e.g.\ \citep[Theorem C.20]{leo17}, with \(\eps = \eta/K\) and \(u(\cdot)\) as the indicator function of \(U\), i.e.\ \(u\mid_U \equiv 1\) and \(u \mid_{\RR^n \setminus U} = 0\).

	Consider \(W: \RR^n \to \RR\) as defined in the previous case. The gradient of \(W\) is uniformly bounded, so there is a fixed point \(y\) such that \(y = x - \tau \overline{\nabla} W(x,y)\). By the same arguments as in case \emph{(ii)}, \({\overline{\nabla} V(x,y) = \overline{\nabla} W(x,y)}\), which implies that \(y\) solves \(y = x - \tau \overline{\nabla} V(x,y)\).
\end{proof}

The third case in \thmref{thm:existence} covers problems where \(V\) is not coercive, which includes linear systems with nonempty kernel and logistic regression problems \citep{lee06} without regularisation.

While the above theorem holds also for the Itoh--Abe methods, there is a much simpler existence result in this case, given in \citep{rii21}, for which continuity of the objective function is sufficient.

\section{Solving the discrete gradient equation}\label{sec:fixed_point}

In the previous section, we proved that the discrete gradient equation~\eqref{eq:dg} admits a solution for \(y\) for all \(\tau > 0\) and \(x \in \RR^n\). In what follows, we propose a relaxed fixed point method for solving~\eqref{eq:dg}, and prove linear convergence rates of the method for the mean value discrete gradient equation. This analysis can trivially be extended to the Itoh--Abe discrete gradient equation, as solving this equation corresponds to solving a succession of mean value discrete gradient equations in one dimension (since all discrete gradients are the same in one dimension, being implicitly defined by the mean value property~\eqref{eq:mean_value}).
\begin{remark}
	We were unable to prove linear convergence for the Gonzalez discrete gradient equation. However, in practice the scheme converged to a solution as quickly as for the mean value case.
\end{remark}

\cite{nor14} show that for sufficiently small time steps, there exists a unique solution to~\eqref{eq:dg} that can be approximated by the fixed point iterations % chktex 2
\begin{equation}\label{eq:fixed_point}
	y^{j+1} = T_\tau(y^j), \quad \mbox {for } j \in \NN \mbox{, where }  \quad T_\tau(y) := x - \tau \overline{\nabla} V(x, y),\; y^0 \in \RR^n
\end{equation}
i.e.\ such that the iterates converge to a fixed point \(y^* = T_\tau(y^*)\), or, equivalently, a solution to~\eqref{eq:dg}. For this, it is assumed that \(\tau\) is less than \(1/(10 L_{\text{DG}})\), where \(L_{\text{DG}}\) is the Lipschitz constant for a given \(x\) of the mapping \(y \mapsto \overline{\nabla} V(x,y)\).

One is often interested in taking larger time steps, and particularly for the optimisation of \hbox{\(L\)-smooth} functions, optimal time steps may be closer to \(2/L\)---see \secref{sec:convergence_rate}. Furthermore, as \thmref{thm:existence} ensures the existence of a solution for arbitrarily large time steps, we would like a constructive method for locating such solutions. We therefore propose to use the relaxed fixed point method, which for \(\theta \in (0,1]\) is given by % chktex 9
\begin{equation}\label{eq:relaxed}
	y^{j+1} = S_{\theta,\tau}(y^j), \qquad  S_{\theta,\tau}(y) := (1-\theta)y + \theta T_\tau(y).
\end{equation}
For \(\theta = 1\), this reduces to~\eqref{eq:fixed_point}. In the remainder of this section, we will prove convergence guarantees of~\eqref{eq:relaxed} for all time steps. In \secref{sec:numerical}, we demonstrate its numerical efficiency.

In the following, we assume that the discrete gradient inherits smoothness and strong convexity properties from the gradient.
\begin{assumption}\label{ass:smooth_convex}
	There are constants \(c_L, c_\mu > 0\) such that
	\begin{enumerate}[(i) ]
		\item \emph{(Smoothness)} If \(V\) is \(L\)-smooth, then for all \(x\in\RR^n\), \(y \mapsto \overline{\nabla} V(x,y)\) is \(c_{L}L\)-Lipschitz continuous.
		\item \emph{(Monotonicity)} If \(V\) is \(\mu\)-convex, then for all \(x, y, z \in \RR^n\), we have
		      \[
			      \inner{\overline{\nabla} V(x, y) - \overline{\nabla} V(x,z), y - z} \geq c_\mu \mu \|y - z\|^2.
		      \]
		      We write \(L_{\text{DG}} := c_{L}L\) and \(\mu_{\text{DG}} := c_\mu \mu\).
	\end{enumerate}
\end{assumption}

It is straightforward to verify these properties for the mean value discrete gradient.
\begin{proposition}\label{prop:mean_value}
	The mean value discrete gradient satisfies \assref{ass:smooth_convex} with \(L_{\text{DG}} = L/2\) and \(\mu_{\text{DG}} =\mu/2\).
\end{proposition}
\begin{proof}
	To show that the first property holds, we write
	\begin{align*}
		\|\overline{\nabla} V(x,y) - \overline{\nabla} V(x,z)\| & \leq \int_0^1\! \|\nabla V(sy + (1-s)x) \!-\! \nabla V(sz + (1-s)x)\| \dif s \\ & \leq L \|y-z\| \int_0^1 \!\! s \dif s = \frac{L}{2} \|y-z\|
	\end{align*}
	Similarly, to show the second property, we write
	\begin{align*}
		\inner{\overline{\nabla}  V(x,y) - \overline{\nabla} & V(x,z) , y - z} \\ &= \int_0^1 \! \frac{1}{s} \biginner{ \nabla V(sy + (1-s)x) - \nabla V(sz + (1-s)x), sy - sz } \dif s
		\\ &\geq \mu \|y-z\|^2  \int_0^1 \! s  \dif s = \frac{\mu}{2} \|y-z\|^2.
	\end{align*}
\end{proof}

To demonstrate that the scheme defined in~\eqref{eq:relaxed} converges to a unique solution of \(y^* = T_\tau(y^*)\), we will use Banach's Fixed Point Theorem (see e.g.\ \citep[Theorem 1.1]{gra03}).
\begin{theorem}[Banach's Fixed Point Theorem]
	Let \((X,d)\) be a non-empty complete metric space. Suppose \(T: X \to X\) is a contraction mapping on \(X\), i.e.\ there is \(q \in [0,1)\) such that % chktex 9
	\[
		d(T(x),T(y)) \leq q d(x,y) \quad \mbox{ for all } x, y \in X.
	\]
	Then \(T\) has a unique fixed point \(y^* \in X\), and for any \(y^0 \in X\), the sequence \({(y^j)}_{j \in \NN}\) defined by \({y^{j+1} = T(y^j)}\) converges linearly to \(y^*\).
\end{theorem}

The following result demonstrates that for convex objective functions, the scheme~\eqref{eq:relaxed} converges to a fixed point \(y^* = T_\tau(y^*)\) for arbitrary time steps \(\tau\).
\begin{theorem}\label{thm:fixed_point}
	If \(V\) is \(L\)-smooth and \(\overline{\nabla}\) satisfies \assref{ass:smooth_convex}, then the mapping \(S_{\theta,\tau}\) defined in~\eqref{eq:relaxed} is a contraction mapping on \(\RR^n\) if either of the following holds.
	\begin{enumerate}[(i) ]
		\item \(\tau < 1/L_{\text{DG}}\) and \(\theta \in [0,1]\).
		\item \(V\) is \(\mu\)-convex and \(\theta \in (0, \min\{1, \frac{2+2\tau \mu_{\text{DG}}}{1 + \tau^2 L_{\text{DG}}^2 + 2 \tau \mu_{\text{DG}}}\})\).
	\end{enumerate}
\end{theorem}

\begin{proof}
	Choose \(y, z \in \RR^n\).
	\emph{Case (i).} We write
	\begin{align*}
		\|S_{\theta,\tau}(z)- S_{\theta,\tau}(y)\| & = \|(1-\theta)(z - y) + \tau \theta (\overline{\nabla} V(x, y) - \overline{\nabla} V(x, z))\| \\ \leq \del{1 - (1 - \tau L_{\text{DG}}) \theta} \|z - y\|.
	\end{align*}
	Thus \(S_{\theta,\tau}\) is a contraction provided that \(\tau < 1/L_{\text{DG}}\).

	\emph{Case (ii).} In a similar fashion, we write
	\begin{align*}
		\|S_{\theta,\tau}(z) & - S_{\theta,\tau}(y)\|^2 = \|(1-\theta)(z - y) + \tau \theta \del{\overline{\nabla} V(x, y) - \overline{\nabla} V(x, z)}\|^2                         \\
		                     & = {(1-\theta)}^2 \|z- y\|^2 +  \tau^2 \theta^2 \|\overline{\nabla} V(x, y) - \overline{\nabla} V(x, z)\|^2                                           \\
		                     & \qquad\qquad\qquad\qquad\qquad - 2 \tau (1-\theta) \theta \inner{z - y, \overline{\nabla} V(x, z) - \overline{\nabla} V(x, y)}                       \\
		                     & \leq \del{ {(1-\theta)}^2 +   \tau^2 \theta^2 L_{\text{DG}}^2 - 2 \tau (1-\theta) \theta \mu_{\text{DG}}} \|z - y\|^2 = \omega(\theta)  \|z - y\|^2.
	\end{align*}
	One can check that the coefficient \(\omega(\theta)\) is less than 1 provided \(\theta\) belongs to the interval stated in the theorem.
	This concludes the proof.
\end{proof}

\begin{remark}
	In the second case of \thmref{thm:fixed_point}, the coefficient \(\omega(\theta)\) is minimised for
	\begin{equation}\label{eq:theta}
		\theta^* =  \frac{1+\tau \mu_{\text{DG}}}{1 + \tau^2 L_{\text{DG}}^2 + 2 \tau \mu_{\text{DG}}} < 1,
	\end{equation}
	which yields the linear convergence rate
	\[
		\|y^{j+1}- y^j\|^2 \leq  \frac{\tau^2(L_{\text{DG}}^2 - \mu_{\text{DG}}^2)}{{(1+\tau \mu_{\text{DG}})}^2 + \tau^2 (L_{\text{DG}}^2 - \mu_{\text{DG}}^2)}   \|y^j - y^{j-1}\|^2.
	\]
	We note from this that the scheme converges faster for smaller time steps and for objective functions with smaller condition numbers \(L/\mu \approx L_{\text{DG}}/ \mu_{\text{DG}} =: \kappa_{\text{DG}}\). Furthermore, if \(\tau = 1/(aL_{\text{DG}})\) for some \(a \geq 1\), where a typical choice is \(a = 1\), then we obtain
	\[
		\theta^* = \frac{1 + \frac{1}{a\kappa_{\text{DG}}}}{1 + \frac{1}{a^2} + \frac{2}{a\kappa_{\text{DG}}}} \geq \frac{a^2}{1+a^2}, \quad \omega(\theta^*) = \frac{1 - \frac{1}{\kappa_{\text{DG}}^2}}{a^2 + \frac{2a}{\kappa_{\text{DG}}}+1} \leq \frac{1}{a^2 + 1}.
	\]
	This shows that the fixed point scheme~\eqref{eq:relaxed} is robust to ill-conditioned problems, both with regards to appropriate choices of \(\theta\) and the rate of convergence.
\end{remark}

In \secref{sec:numerical_fixed_point}, we compare the efficiency of the above scheme for different \(\theta\) and of the built-in solver \texttt{scipy.optimize.fsolve} in Python.

\section{Analysis of time steps for discrete gradient methods}\label{sec:analysis_time_step}

In this section, we study the implicit dependence of \(x^{k+1}(\tau)\) on the choice of time step \(\tau\). We first establish a uniqueness result for the mean value and Itoh--Abe discrete gradient methods. Then we restrict our focus to Itoh--Abe methods, where we ascertain bounds on optimal time steps with respect to the decrease in \(V\), for \(L\)-smooth, convex functions as well as strongly convex functions.

\subsection{Uniqueness}

\begin{lemma}
	If \(V\) is \(\mu\)-convex, then the solution \(y\) to the discrete gradient equation~\eqref{eq:dg} is unique for the mean value discrete gradient and the Itoh--Abe discrete gradient.
\end{lemma}

\begin{proof}
	We first consider the mean value discrete gradient. Suppose \(y_1\), \(y_2\), solve \(y_i = x - \tau \overline{\nabla} V(x,y_{i})\), \(i = 1,2\). Then it follows from \propref{prop:mean_value} that
	\begin{align*}
		\|y_1- y_2\|^2 & = \tau \biginner{\frac{x-y_2}{\tau} - \frac{x-y_1}{\tau}, y_1 - y_2} = \tau \inner{\overline{\nabla} V(x, y_2) - \overline{\nabla} V(x,y_1), y_1 - y_2} \leq 0.
	\end{align*}

	Furthermore, as the Itoh--Abe discrete gradient method is a succession of scalar updates, each corresponding to a 1D mean value discrete gradient update, uniqueness follows.
\end{proof}

\subsection{Implicit dependence on the time step for Itoh--Abe methods}

For the remainder of the section, we restrict our focus to Itoh--Abe methods.  We fix a starting point \(x\), direction \(d \in S^{n-1}\) and time step \(\tau\), and study the solution \(y\) to
\begin{equation}\label{eq:dg_scalar}
	y = x - \alpha d, \qquad \mbox{where }\alpha \neq 0 \mbox{ solves } \quad \alpha =  - \tau \frac{V(x - \alpha d) - V(x)}{\alpha}.
\end{equation}
By the analysis in \secref{sec:existence}, there exists a solution \(y\) for all \(\tau > 0\). For convenience and to exclude the case \(y  = x\), we assume \(\inner{\nabla V(x) , d} > 0\). For notational brevity, we rewrite the optimisation problem in terms of a scalar function \(f\), i.e.\ solve
\begin{equation}\label{eq:dg_scalar2}
	\dfrac{f(\alpha)}{\alpha^2} = -\frac{1}{\tau}, \quad \mbox{ where } f(\alpha) := V(x - \alpha d) - V(x).
\end{equation}

For optimisation schemes with a time step \(\tau\), it is common to assume that the distance between \(x\) and \(y\) increases with the time step. For explicit schemes, this naturally holds. However, for implicit schemes, such as the discrete gradient method, this is not always the case. We demonstrate this with a simple example in one dimension.
\begin{example}
	Define \(V(x) := -x^3\) and \(x = 0\). For all \(\tau > 0\), \eqref{eq:dg_scalar} is solved by \(y = \frac{1}{\tau}\). Then, as \(\tau \to 0\), we have \(y\to\infty\), and as \(\tau \to \infty\), we have \(y \to x\). % chktex 2
\end{example}
The above example illustrates that for nonconvex functions, decreasing the time step might lead to a larger step \(y \mapsfrom x\) and vice versa.

We now show that for convex functions, the distance \(\|y - x\|\) does increase with \(\tau\). Set
\[
	R = \sup\set{r \; : \; V(x-\alpha d) < V(x) \mbox{ for all } \alpha \in (0, r)}.
\]
By the assumption that \(\inner{\nabla V(x), d} > 0\), we have \(R > 0\).

\begin{proposition}\label{prop:convex_bijection}
	If \(V\) is convex, then there is a well-defined, continuous, and strictly increasing mapping \(\tau \mapsto \alpha(\tau)\), such that \(\alpha(\tau)\) solves~\eqref{eq:dg_scalar2} for \(\tau\). Furthermore, the mapping is bijective from \((0,\infty)\) to \((0,R)\).
\end{proposition}
\begin{proof}
	We first establish continuity and strict monotonicity. We use the alternative characterisation of convex functions in one dimension, which states that
	\[
		\alpha \mapsto  \dfrac{f(\alpha)- f(0)}{\alpha} = \dfrac{f(\alpha)}{\alpha}
	\]
	is monotonically nondecreasing in \(\alpha\). Since \(f(\alpha) < 0\) for \(\alpha \in (0, R)\), it follows that \(\alpha \mapsto f(\alpha)/\alpha^2\) is strictly increasing on \((0, R)\). We can thus apply e.g.\ the implicit function theorem for strictly monotone functions \citep[Theorem 1H.3]{don14} to conclude that the mapping \(\tau \mapsto \alpha(\tau)\) is continuous for all \(\tau > 0\).

	Furthermore, for any \(\tau_2 > \tau_1>0\) and corresponding solutions \(\alpha_1\), \(\alpha_2\), using~\eqref{eq:dg_scalar2}, we get
	\[
		\frac{f(\alpha_2)}{\alpha_2^2} = -\frac{1}{\tau_2} > -\frac{1}{\tau_1} = \frac{f(\alpha_1)}{\alpha_1^2}
	\]
	Since \(\alpha \mapsto f(\alpha)/\alpha^2\) is strictly increasing, it follows that \(\alpha_2 > \alpha_1\). Hence \(\tau \mapsto \alpha(\tau)\) is strictly increasing.

	Next, we show that \(\alpha(\tau) \to 0\) as \(\tau \to 0\). This can be seen by inspecting
	\[
		\dfrac{f(\alpha(\tau))}{\alpha(\tau)} = -\dfrac{\alpha(\tau)}{\tau}.
	\]
	The left-hand side is bounded by the derivative \(f'(0) = - \inner{\nabla V(x), d}\). Hence, as \(\tau\) goes to zero, \(\alpha(\tau)\) must also go to zero to prevent the right-hand side from blowing up.

	Last, we show that \(\alpha(\tau) \to R\) as \(\tau \to \infty\). By inspecting
	\[
		\dfrac{f(\alpha(\tau))}{\alpha{(\tau)}^2} = -\dfrac{1}{\tau},
	\]
	we see that as \(\tau \to \infty\), the right-hand side goes to zero, so either \(f(\alpha(\tau)) \to 0\) or \(\alpha{(\tau)}^2 \to \infty\). There are two cases to consider for \(R\). If \(R < \infty\),
	then by convexity of \(f\) and the definition of \(R\), \(f(\alpha) < 0\) for \(\alpha \in (0,R)\) and \(f(\alpha) > 0\) for \(\alpha > R\), so \(\alpha(\tau) \to R\). Otherwise, if \(R = \infty\), then there exists an \(\eps > 0\) such that \(f(\alpha) < -\eps\) for all \(\alpha > 0\), from which it follows that \(\alpha{(\tau)}^2 \to \infty = R\). This concludes the proof.
\end{proof}

\begin{remark}
	The above proposition can also be shown to hold for non-differentiable, convex functions, by replacing the derivative with a subgradient.
\end{remark}

\subsection{Lipschitz continuous gradients}

The remainder of this section is devoted to deriving bounds on optimal time steps, with respect to the decrease in the objective function when the objective function is \(L\)-smooth or \(\mu\)-convex. We first consider \(L\)-smooth functions, and show that any time step \(\tau < 2/L\) is suboptimal. We recall the scalar function \(f(\alpha) = V(x - \alpha d) - V(x)\).
\begin{lemma}\label{lem:smooth}
	If \(V\) is convex and \(L\)-smooth, then \(\tau \mapsto f(\alpha(\tau))\) is decreasing for \(\tau < 2/L\).
\end{lemma}

\begin{proof}
	Suppose \(\alpha\) solves~\eqref{eq:dg_scalar2} for \(\tau < 2/L\). Let \(\lambda \in (\tau L / 2, 1)\), and plug in \(0\) and \(\alpha/\lambda\) for \(y\) and \(x\) respectively in \propref{prop:smooth} \emph{(ii)} to get, after rearranging,
	\[
		\lambda f(\alpha/\lambda)  \leq  f(\alpha) + \dfrac{(1-\lambda) L}{2 \lambda} \alpha^2.
	\]
	Plugging in~\eqref{eq:dg_scalar2}, we get
	\[
		f(\alpha/\lambda)  \leq  \del{\dfrac{1}{\lambda} - \dfrac{(1-\lambda)\tau L}{2 \lambda^2}  } f(\alpha).
	\]
	We  show that \(f(\alpha/\lambda) < f(\alpha)\), i.e.\ that \(\lambda - (1-\lambda)\tau L/2 > \lambda^2\). By solving the quadratic expression, we find this holds when \(\lambda \in (\tau L/2, 1)\). Thus \(f(\alpha/\lambda) < f(\alpha)\), and, as \(f\) is convex, it follows that \(f\) is decreasing on \([\alpha, \alpha/\lambda]\). By \propref{prop:convex_bijection}, \(\tau \mapsto f(\alpha(\tau))\) is therefore decreasing on \((0,2/L)\).
\end{proof}

\subsection{Strong convexity}

We now assume \(\mu > 0\) and show that for strongly convex functions, any time step \(\tau > 2/\mu\) yields a suboptimal decrease.

\begin{lemma}\label{lem:convex}
	If \(V\) is \(\mu\)-convex, where \(\mu > 0\), then \(\tau \mapsto f(\alpha(\tau))\) is strictly increasing for \(\tau > 2/\mu\).
\end{lemma}

\begin{proof}
	Let \(\alpha\) solve~\eqref{eq:dg_scalar2} for \(\tau > 2/\mu\). Fix \(\lambda \in (2/(\tau \mu), 1)\), and plug in \(0\) and \(\alpha\) for \(y\) and \(x\) respectively in \defnref{defn:strongly_convex} \emph{(ii)} to get, after rearranging,
	\[
		f(\lambda\alpha)  \leq  \lambda f(\alpha)  - \dfrac{\mu \lambda (1-\lambda)}{2} \alpha^2.
	\]
	Plugging in~\eqref{eq:dg_scalar2} gives us
	\[
		f(\lambda\alpha) \leq \del{  \lambda  + \dfrac{\tau \mu \lambda (1-\lambda)}{2}} f(\alpha).
	\]
	We want to show that \(f(\lambda \alpha) < f(\alpha)\), i.e.\ that \(\lambda  + \tau \mu \lambda (1-\lambda)/2 > 1\). By rearranging and solving the quadratic expression, we find that this is satisfied if \(\lambda \in (2/(\tau \mu), 1)\). The result follows from convexity of \(f\) and \propref{prop:convex_bijection}.
\end{proof}
\begin{remark}
	This result also holds for strongly convex, non-differentiable functions.
\end{remark}

\section{Convergence rate analysis}\label{sec:convergence_rate}

In this section we derive convergence rates for \(L\)-smooth functions, \(\mu\)-convex functions, and functions that satisfy the Polyak--{\L}ojasiewicz (P{\L}) inequality. We follow the arguments in \citep{bec13, nes12}, on convergence rates of coordinate descent. We assume throughout that the iterates \({(x^{k})}_{k\in\NN}\) solve the discrete gradient equation exactly, and leave the impact of inexact updates on the convergence rates for future work.

We recall the notation in~\eqref{eq:expectation}, \(V_{k} := \EE [V(x^{k})]\), where \(V_{k} = V(x^{k})\) for deterministic methods. Estimates of the following form will be crucial to the analysis.
\begin{equation}\label{eq:estimate_main}
	\beta \del{V(x^{k}) - \EE_{d^{k+1}} [V(x^{k+1})]} \geq \|\nabla V(x^{k})\|^2.
\end{equation}
We first provide this estimate for each of the four methods. We assume throughout that the time steps \({(\tau_{k})}_{k \in \NN}\) satisfy arbitrary bounds~\eqref{eq:time_step}.

We consider coordinate-wise Lipschitz constants for the gradient of \(V\) as well as a directional Lipschitz constant. For \(i =1, \ldots, n\), we suppose \(\partial_i V : \RR^n \to \RR^n\) is Lipschitz continuous with Lipschitz constant \(\overline{L}_i \leq L\). We denote by \(\overline{L}_{\Sum}\) the \(\ell^2\)-norm of the coordinate-wise Lipschitz constants, \({\overline{L}_{\Sum} = \sqrt{\sum_{i=1}^n  \overline{L}_i^2} \in [L , \sqrt{n} L]}\).

Furthermore, for a direction \(d \in S^{n-1}\), we consider the Lipschitz constant \(L_d \leq L\), such that for all \(x \in \RR^n\) and \(\alpha \in \RR\), we have
\[
	|\inner{\nabla V(x + \alpha d), d} - \inner{\nabla V(x), d}| \leq L_d |\alpha|.
\]
For the Itoh--Abe discrete gradient method or when \(\Xi\) only draws from the standard coordinates, we write \(L_i\) instead of \(L_{e^{i}}\). We define \(L_{\max} \leq L\) to be the supremum of \(L_d\) over all \(d\) in the support of the probability density function of \(\Xi\). That is, \(L_{\max} \geq L_d\) for all \(d \sim \Xi\). In the case when \(\Xi\) draws from a restricted set, such as the standard coordinates, \(L_{\max}\) can be notably smaller than \(L\). Hereby, we refine the \(L\)-smoothness property in \propref{prop:smooth} \emph{(i)} to
\begin{equation}\label{eq:lipschitz_direction}
	V(x + \alpha d) - V(x) \leq \alpha \inner{\nabla V(x), d} + \dfrac{L_d}{2} \alpha^2 \leq \alpha \inner{\nabla V(x), d} + \dfrac{L_{\max}}{2} \alpha^2,
\end{equation}
for all \(\alpha \in \RR\) and \(d\) in the support of the density of \(\Xi\) \citep[Lemma 3.2]{bec13}.

\begin{table}
	\caption{Estimates of \(\beta\), as well as optimal time steps \(\tau^*\) and \(\beta^*\), with \(\zeta\) given in~\eqref{eq:angle_bound}.\label{tab:time_step}}
	\begin{tabular}{@{}cccc@{}}
		\toprule
		Discrete gradient method & \(\beta\)                                                       & \(\tau^*\)                & \(\beta^*\)               \\
		\midrule
		Gonzalez                 & \(\displaystyle 2 (1/\tau_k + L^2 \tau_k/2)\)                   & \(\sqrt{2}/L\)            & \(2 \sqrt{2}L \)          \\
		Mean value               & \(\displaystyle 2 ( 1/\tau_k + L^2 \tau_k/4 )\)                 & \(2/L\)                   & \(2L\)                    \\
		Itoh--Abe                & \(\displaystyle 2 (1/\tau_k + \overline{L}_{\Sum}^2 \tau_{k})\) & \(1/\overline{L}_{\Sum}\) & \(4 \overline{L}_{\Sum}\) \\
		Randomised Itoh--Abe     & \(\displaystyle \tau_k {(1/\tau_k + L_{\max}/2)}^2/\zeta\)      & \(2/L_{\max}\)            & \(2 L_{\max}/\zeta\)      \\
		\bottomrule
	\end{tabular}
\end{table}

\begin{lemma}\label{lem:estimate_main}
	If \(V\) is \(L\)-smooth, then the three discrete gradient methods and the randomised Itoh--Abe method satisfy~\eqref{eq:estimate_main} with values for \(\beta\) given in \tabref{tab:time_step}.
\end{lemma}

A proof of this lemma is given in the \appref{sec:estimate_main}. Note that these estimates do not require convexity of \(V\).

\subsection{Optimal time steps and estimates of \texorpdfstring{$\beta$}{\textbeta}}\label{sec:optimal_time_step} % chktex 46

Lower values for \(\beta\) in~\eqref{eq:estimate_main} correspond to better convergence rates, as can be seen in Theorems~\ref{thm:rate1} and~\ref{thm:rate2}. In what follows, we briefly discuss the time steps that yield minimal values of \(\beta\), denoted by \(\tau^*\) and \(\beta^*\) in~\tabref{tab:time_step}.

For the Gonzalez and mean value discrete gradient methods, it is natural to compare rates to those of explicit gradient descent, which has the estimate \(\beta^* = 2L\)  \citep[Section 2.1.5]{nes04}. Hence, the mean value discrete gradient method recovers the optimal rates of gradient descent, while the estimate for the Gonzalez discrete gradient is worse by a factor of \(\sqrt{2}\).

In contrast, for the proximal point method, it can be shown that \(\beta = 2/\tau\) and can thus be made arbitrarily small by taking larger time steps. This follows from the fact that the proximal point method with time step \(\tau\) applied to \(V\) is equivalent to explicit gradient descent with time step \(\tau\) applied to the Moreau--Yosida regularisation of \(V\) with parameter \(\tau\), \(V_\tau\), which is \hbox{\(1/\tau\)-smooth} \citep[Section 4.2]{cha16}. It follows that \(\beta = 2/\tau\).

For the Itoh--Abe discrete gradient method, we compare its rates to those obtained for cyclic coordinate descent (CD) schemes in \citep[Theorem 3]{wri15} and \citep[Lemma 3.3]{bec13}, \(\beta^* = 8 \sqrt{n} L\), where we have set their parameters \(L_{\max}\) and \(L_{\min}\) to \(\sqrt{n} L\). Hence, the estimate for the Itoh--Abe discrete gradient method is stronger, being at most half that of CD, even in the worst-case scenario \(\overline{L}_{\Sum} = \sqrt{n} L\).
\begin{remark}
	Note however that we can improve the estimate for the CD scheme to recover the same rate. See \appref{sec:CD}.
\end{remark}

We give one motivating example for considering the parameter \(\overline{L}_{\Sum}\).
\begin{example}
	Let \(V\) be a least squares problem \(V(x) = \|Ax - f\|^2/2\). We then have
	\begin{equation}\label{eq:rank}
		\overline{L}_{\Sum} \leq \sqrt{\rank(A)} L.
	\end{equation}
	Thus, for low-rank system where \(\rank(A) \ll  n\), the convergence speed of the Itoh--Abe discrete gradient method improves considerably.

	To derive~\eqref{eq:rank}, one can show that \(L = \|A^*A\|_2\) and \(\overline{L}_{\Sum} = \|A^* A\|_F\), where \(\|\cdot\|_2\) and \(\|\cdot\|_F\) denote the operator norm and Frobenius norm. The bound follows from the properties  \(\|B\|_F \leq \sqrt{\rank(B)} \|B\|_2\) \citep[Table 6.2]{hig02} and \(\rank(A^*A) = \rank(A)\) \citep[Statement 4.5.4]{mey00}.
\end{example}

We compare the rates for the randomised Itoh--Abe methods to randomised coordinate descent (RCD). Recall that when \(\Xi\) is the uniform distribution on the coordinates \({(e^{i})}_{i=1}^n\) or on the unit sphere \(S^{n-1}\), we have \(\zeta = 1/n\). This gives us \(\beta^* = 2 n L_{\max}\) for the randomised Itoh--Abe methods, which is the optimal bound for randomised coordinate descent \citep{wri15}.

\subsection{Lipschitz continuous gradients}

For the next result, we use the notation \(R(x^0) := \diam \{x \in \RR^n \; : \; V(x) \leq V(x^0)\}\). \(R(x^0)\) is bounded, provided \(V\) is coercive.

\begin{theorem}\label{thm:rate1}
	Let \(V\) be an \(L\)-smooth, convex, coercive function. Then for all four methods, \(\beta\) given in \tabref{tab:time_step}, and \(V^* := \min_x V(x)\), we have
	\begin{align}\label{eq:1/k}
		V_k - V^* \leq \dfrac{\beta R{(x^0)}^2 }{k + 2\frac{\beta}{L}}.
	\end{align}
\end{theorem}
\begin{proof}
	Let \(x^*\) be a minimiser of \(V\). By respectively convexity, the Cauchy-Schwarz inequality, and \lemref{lem:estimate_main}, we have
	\begin{align*}
		{(V(x^{k}) - V^*)}^2 & \leq \biginner{\nabla V(x^{k}), x^{k} - x^{*}}^2  \leq  \| \nabla V(x^{k}) \|^2 \|x^{k} - x^{*} \|^2 \\ & \leq  \beta R{(x^0)}^2 (V(x^{k}) - \EE_{d^{k+1}}[V(x^{k+1})]).
	\end{align*}
	Taking expectation on both sides with respect to \(\xi_{k}\), and applying Jensen's inequality \citep[Theorem 3.3]{rud87}, we obtain
	\[
		{(V_k - V^*)}^2 \leq  \beta R{(x^0)}^2 (V_k - V_{k+1}).
	\]
	By monotonicity of \(V_k\) and following the steps in the proof of \citep[Theorem 1]{nes12}, we obtain
	\begin{align*}
		V_k - V^* \leq \dfrac{\beta  R{(x^0)}^2 }{k + \beta  \frac{R{(x^0)}^2}{V(x^0) - V^*}}.
	\end{align*}
	To eliminate dependence on the starting point, we derive \(V(x^0) - V^* \leq  \tfrac{L}{2}R{(x^0)}^2\) from \propref{prop:smooth} \emph{(i)},
	which gives us the desired result~\eqref{eq:1/k}.
\end{proof}

\subsection{The Polyak--{\L}ojasiewicz inequality}

The next result states that for \(L\)-smooth functions that satisfy the P{\L} inequality, we achieve a linear convergence rate. A function is said to satisfy the P{\L} inequality with parameter \(\mu > 0\) if, for all \(x \in \RR^n\),
\begin{equation}\label{eq:PL}
	\frac{1}{2} \|\nabla V(x) \|^2 \geq \mu \del{V(x) - V^*}.
\end{equation}
Originally formulated by~\cite{pol63}, it was recently shown that this inequality is weaker than other properties commonly used to prove linear convergence \citep{bol10, bol17, csi17, kar16, nec18}. This is useful for extending linear convergence rates to functions that are not strongly convex, including some nonconvex functions.

\begin{proposition}[\cite{kar16}]
	Let \(V\) be \(\mu\)-convex. Then \(V\) satisfies the P{\L} inequality~\eqref{eq:PL} with parameter \(\mu\).
\end{proposition}

We now proceed to the main result of this subsection.
\begin{theorem}\label{thm:rate2}
	Let \(V\) be \(L\)-smooth and satisfy the P{\L} inequality~\eqref{eq:PL} with parameter \(\mu\). Then, for \(\beta\) given in \tabref{tab:time_step}, the three discrete gradient methods and the randomised Itoh--Abe method obtain the linear convergence rate
	\begin{equation}\label{eq:linear_rate}
		V_k - V^* \leq \del{1 - \dfrac{2 \mu}{\beta} }^k\del{V(x^0) - V^*}.
	\end{equation}
\end{theorem}

\begin{proof}
	This is a standard argument, see e.g.\ \citep{kar16}. Combining the P{\L} inequality~\eqref{eq:PL} with the estimate in \lemref{lem:estimate_main}, and taking expectation with respect to \(\xi^{k}\) on both sides, we get
	\begin{align*}
		V(x^{k}) - \EE_{d^{k+1}} [V(x^{k+1})] \geq \dfrac{2 \mu}{\beta}(V(x^{k}) - V^*)  \implies  V_{k+1} - V^* \leq \del{1 - \dfrac{2 \mu}{\beta}} (V_k - V^*).
	\end{align*}
\end{proof}

\section{Preconditioned discrete gradient method}\label{sec:preconditioned}

We briefly discuss the generalisation of the discrete gradient method~\eqref{eq:dg_method} to a preconditioned version
\begin{equation}\label{eq:dg_preconditioned}
	x^{k+1} = x^k - A_k \overline{\nabla} V(x^k, x^{k+1}),
\end{equation}
where \({(A_{k})}_{k \in \NN} \subset \RR^{n \times n}\) is a sequence of positive-definite matrices. Denoting by \(\sigma_{1, k}\) and \(\sigma_{n,k}\) the smallest and largest singular values of \(A_k\) respectively, we have, for all \(x\),
\[
	\sigma_{1, k} \|x\| \leq \|A_k x\| \leq \sigma_{n, k} \|x\|.
\]
It is straightforward to extend the results in \secref{sec:existence} and \secref{sec:convergence_rate} to this setting, under the assumption that there are \(\sigma_{\max} \geq \sigma_{\min} > 0\) such that \(\sigma_{\min} \leq \sigma_{1, k}, \sigma_{n,k} \leq \sigma_{\max}\) for all \(k \in \NN\).

There are several motivations to precondition. In the context of geometric integration, it is typical to group the gradient flow system~\eqref{eq:gradient_flow} with the more general dissipative system
\[
	\dot{x} = - A(x) \nabla V(x),
\]
where \(A(x) \in \RR^{n \times n}\) is positive-definite for all \(x \in \RR^n\) \citep{qui96}. This yields numerical schemes of the form~\eqref{eq:dg_preconditioned}, where we absorb \(\tau_k\) into \(A_k\). There are optimisation problems in which the time step \(\tau_k\) should vary for each coordinate. This is, for example, the case when one derives the SOR method from the Itoh--Abe discrete gradient method \citep{miy18}.

More generally, if one has coordinate-wise Lipschitz constants for the gradient of the objective function, it may be beneficial to scale the coordinate-wise time steps accordingly.

\section{Numerical experiments}\label{sec:numerical}

In this section, we apply the discrete gradient methods to various test problems. The code for the figures has been implemented in Python and MATLAB and is available at \url{https://github.com/riis-academic/discrete-gradient-method-smooth-optimisation/}. For solving the discrete gradient equation~\eqref{eq:dg_method} with the Gonzalez and mean value discrete gradients, we use the fixed point method~\eqref{eq:relaxed} detailed in \secref{sec:fixed_point} and tested numerically in \secref{sec:numerical_fixed_point} under the label `\textbf{R}'. For solving~\eqref{eq:dg_method} for the Itoh--Abe method, we use the built-in solver \texttt{scipy.optimize.fsolve} in Python.  % chktex 12

For these test problems, we generally assume knowledge of the Lipschitz constant \(L\) when setting time steps. However, in \secref{sec:TV} we assess Itoh--Abe methods with a wider range of time steps.

\subsection{Setup}

We use the following time steps for the different methods, unless otherwise specified. For the mean value discrete gradient method, we use \(\tau_{\text{MV}} = 2/L\), for the Gonzalez discrete gradient method, we use \(\tau_{\text{G}} = 2/L\), and for the Itoh--Abe methods, we use the coordinate-dependent time steps \({\tau_{\text{IA},i} = \tau_{\text{RIA},i} = 2/L_i}\). Note that the time steps for the Itoh--Abe discrete gradient method are not the optimal choice suggested in \tabref{tab:time_step}, but were heuristically optimal for the test problems we considered. An analysis of coordinate-dependent time step strategies is an open topic for future research.

In figure captions and legends, the abbreviations \emph{CIA} and \emph{RIA} refer respectively to the (cyclic) Itoh--Abe discrete gradient method and the randomised Itoh--Abe method drawing uniformly from the standard coordinates. For the sake of comparison, we define one iterate of the randomised Itoh--Abe methods as \(n\) scalar updates, so that the computational time is comparable to the standard Itoh--Abe discrete gradient method.

Unless otherwise specified, matrices and vectors are constructed by standard Gaussian draws. To provide the matrix with a given condition number, we compute its singular value decomposition and linearly transform its eigenvalues accordingly.

\subsection{Linear systems}

We first solve linear systems of the form
\begin{equation}\label{eq:linear_problem}
	\min_{x \in \RR^n} V(x) = \frac{1}{2} \|Ax-b\|^2, \quad A \in \RR^{n\times n},\quad b \in \RR^n.
\end{equation}
For linear systems, the Gonzalez and the mean value discrete gradient are both given by \({\overline{\nabla} V(x,y) = \nabla V(\tfrac{x+y}{2}) = A^*(A\tfrac{x+y}{2} - b)}\), so we consider these jointly. As discussed previously, the Itoh--Abe methods reduce to SOR methods for linear systems and are therefore explicit.

\subsubsection{Effect of the condition number}

\begin{figure}[ht]
	\centering
	\includegraphics[height=4cm]{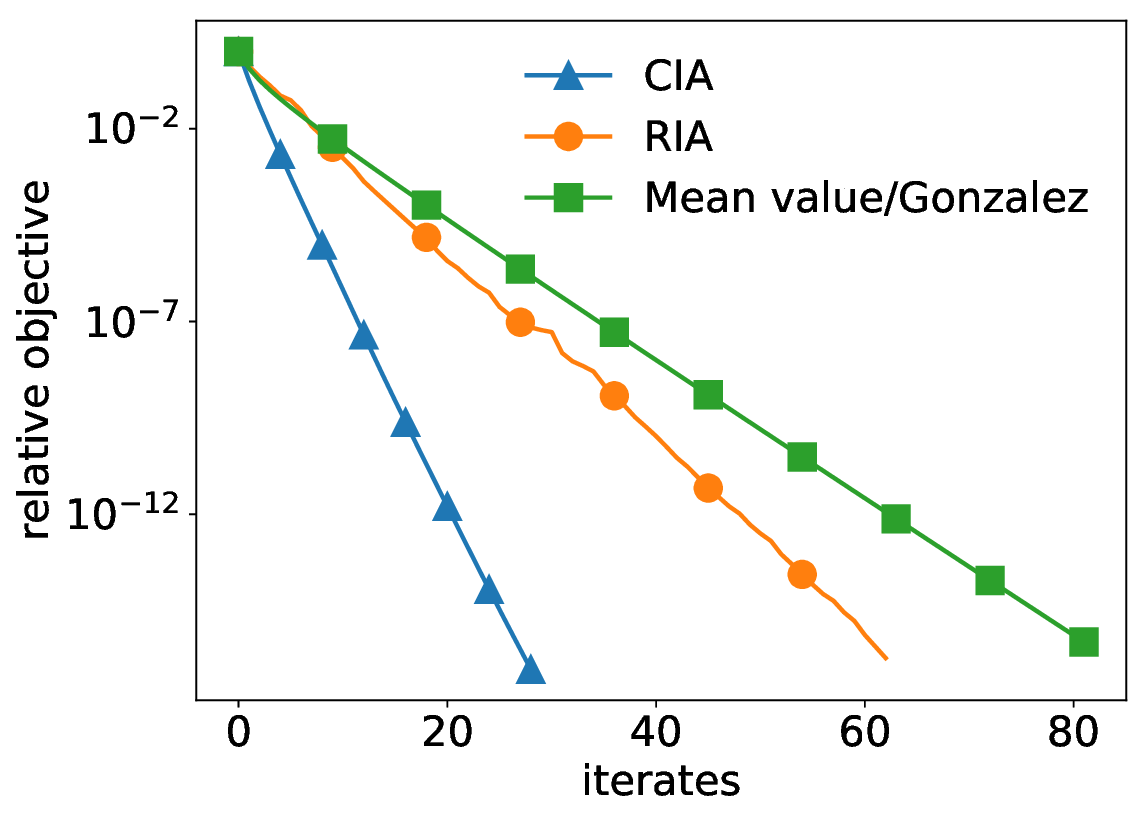}
	\includegraphics[height=4cm]{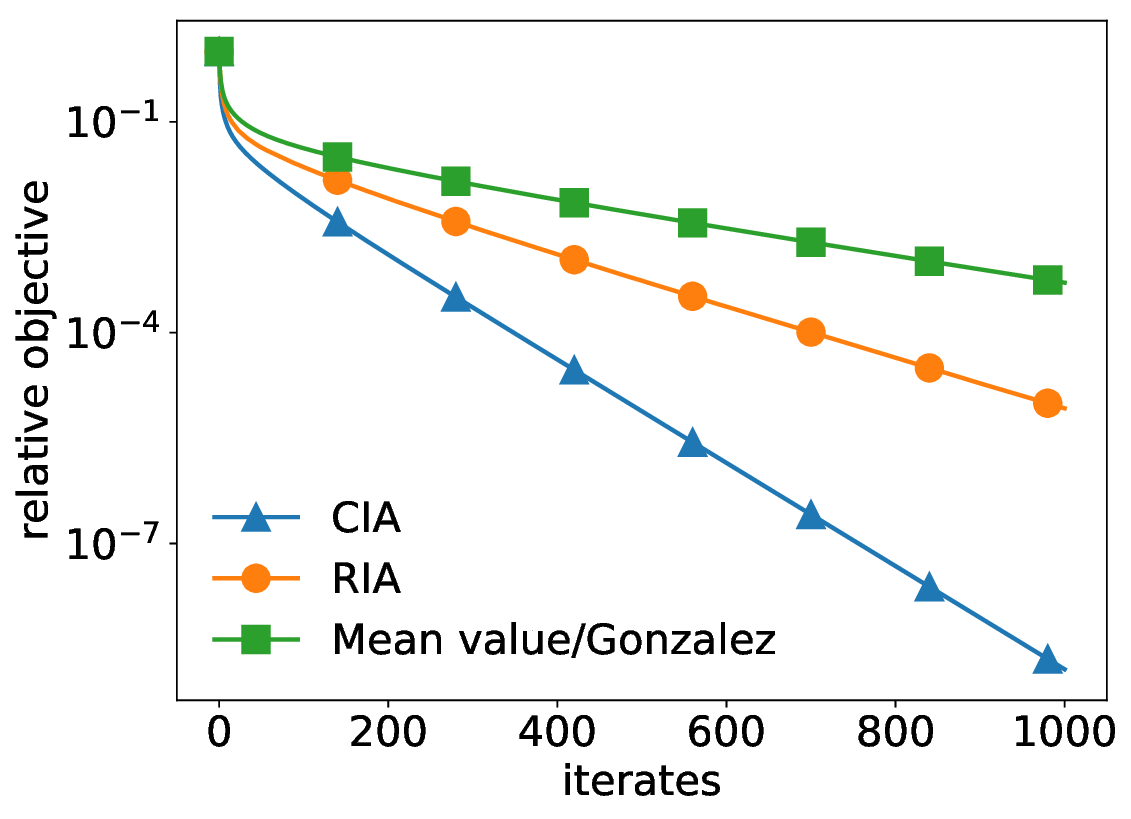}
	\caption{DG methods for linear systems with condition number \(\kappa = 10\) (\textbf{left}) and \(\kappa = 1,000\) (\textbf{right}). Convergence rate plotted as relative objective \([V(x^{k}) - V^*]/[V(x^0) - V^*]\). Linear rate is observed for all methods and is sensitive to condition number.\label{fig:linear}}
\end{figure}

We set \(n = 500\) and consider two linear systems respectively with a low condition number \(\kappa = L/\mu = 10\) and a high one \(\kappa = 1,000\). In both cases, we set \(x^0 = 0\). See \figref{fig:linear} for the results for both cases.

\subsubsection{Sharpness of rates}

\begin{figure}[ht]
	\centering
	\includegraphics[height=4cm]{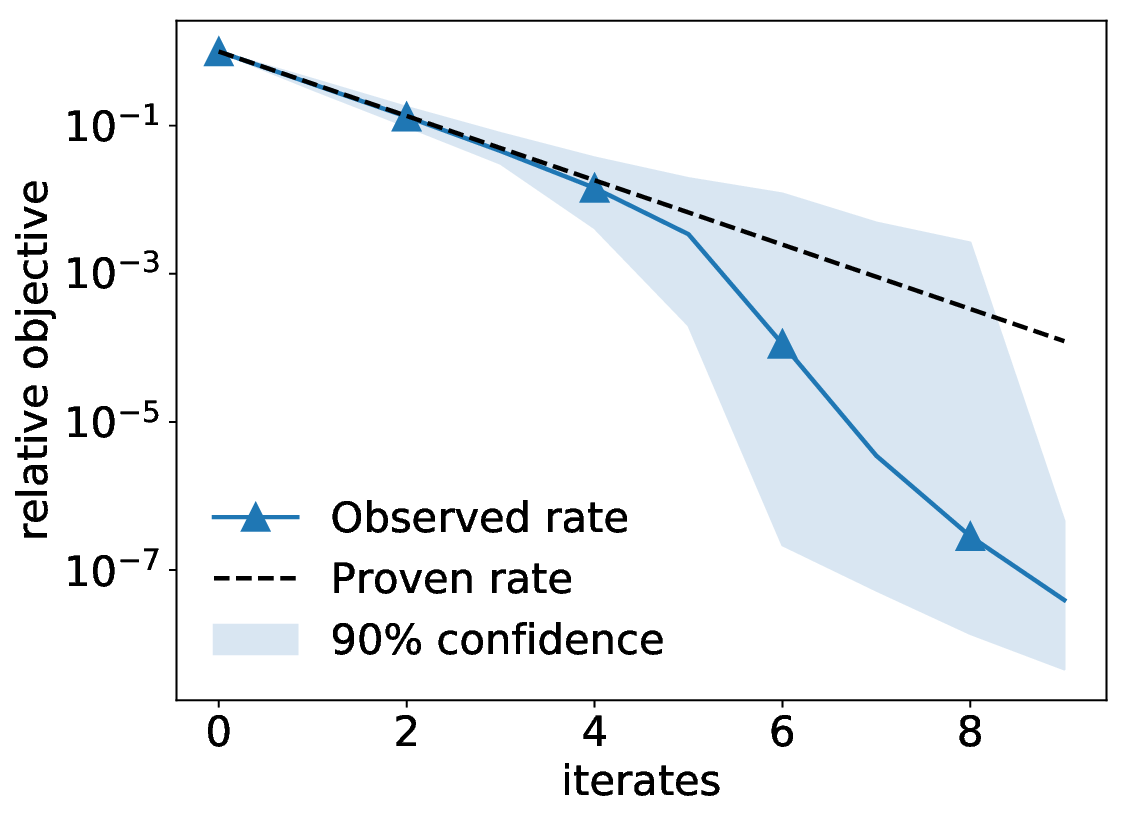}
	\includegraphics[height=4cm]{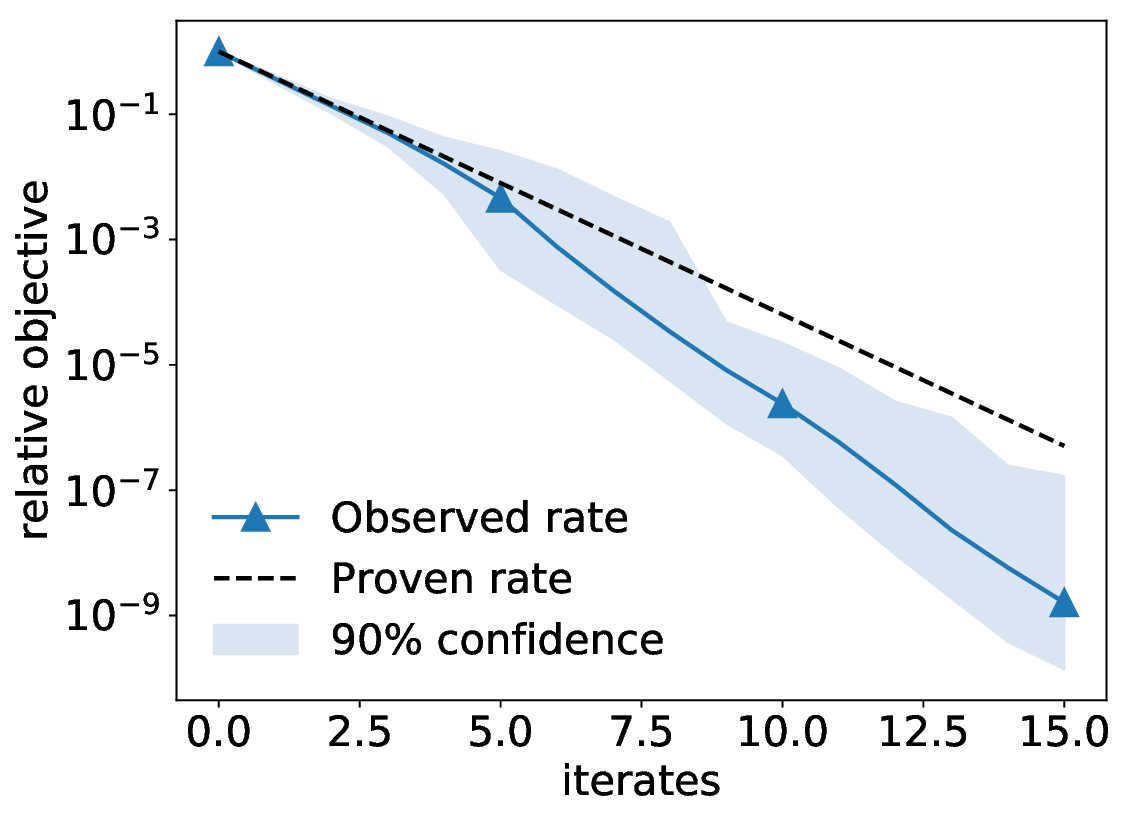}
	\caption{Comparison of observed convergence rate with theoretical convergence rate~\eqref{eq:linear_rate}, for randomised Itoh--Abe method applied to linear system with condition numbers \(\kappa = 1.2\) (\textbf{left}) and \(\kappa = 10\) (\textbf{right}). Average convergence rate and confidence intervals as estimated from 100 runs on the same system. The sharpness of the proven convergence rate is observed in both cases.\label{fig:stats}}

\end{figure}

We test the sharpness of the convergence rate~\eqref{eq:linear_rate} for the randomised Itoh--Abe method. To do so, we run 100 instances of the numerical experiment in the previous subsection and plot the mean convergence rate and 90\(\%\)-confidence intervals, and compare the results to the proven rate. We do this for two condition numbers, \(\kappa = 1.2\) and \(10\). The results are presented in \figref{fig:stats}. These plots suggest that the proven convergence rate estimate is sharp for the randomised Itoh--Abe method.

\subsubsection{Linear system with kernel}

\begin{figure}[ht]
	\includegraphics[height=4cm]{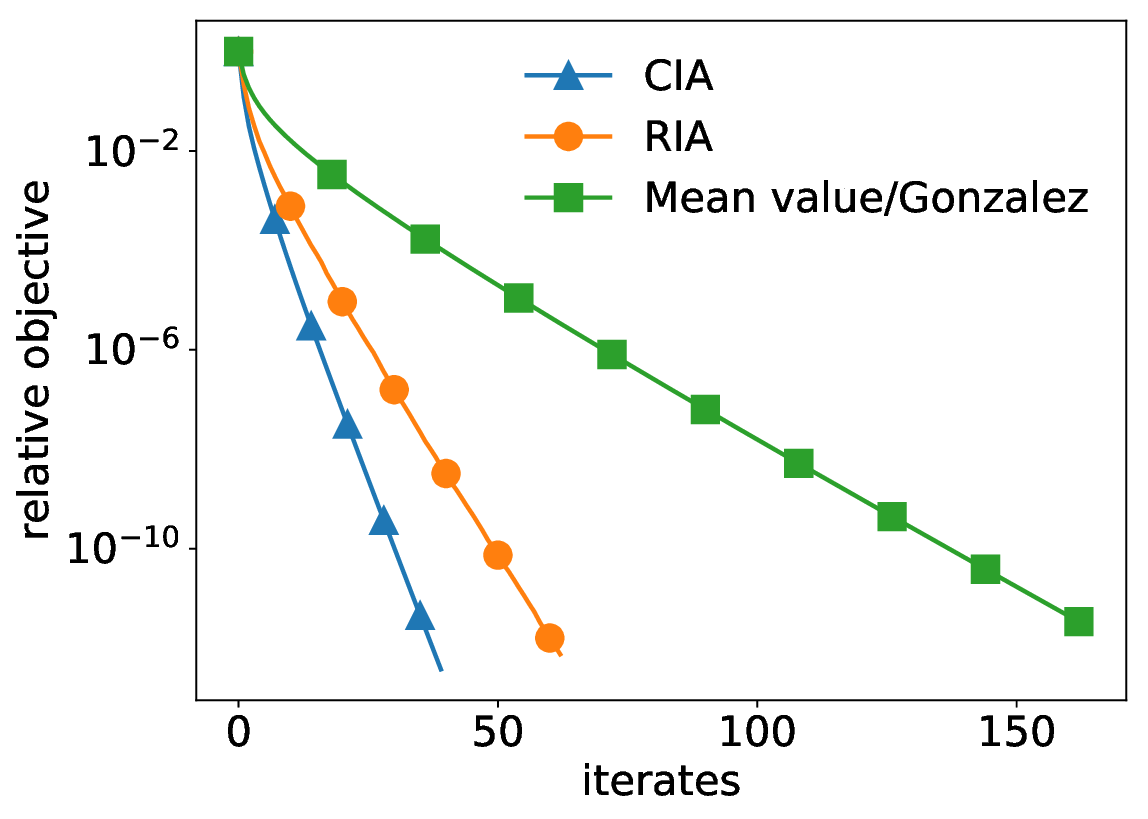}
	\caption{DG methods for linear systems with nontrivial kernel, and convergence rate plotted as relative objective. The function is not strongly convex but satisfies the P{\L} inequality, yielding linear convergence rates.\label{fig:kernel}}
\end{figure}

Next we consider linear systems where the operator \(A\) has a nontrivial kernel, meaning that the objective function is not strongly convex, but nevertheless satisfies the P{\L} inequality. We let \(A \in \RR^{m \times n}\) and \(b \in \RR^m\), where \(n = 800\) and \(m = 400\), meaning the kernel of \(A\) has dimension 400.  See \figref{fig:kernel} for results.

\subsubsection{A note of caution}

\begin{figure}[ht]
	\includegraphics[height=4cm]{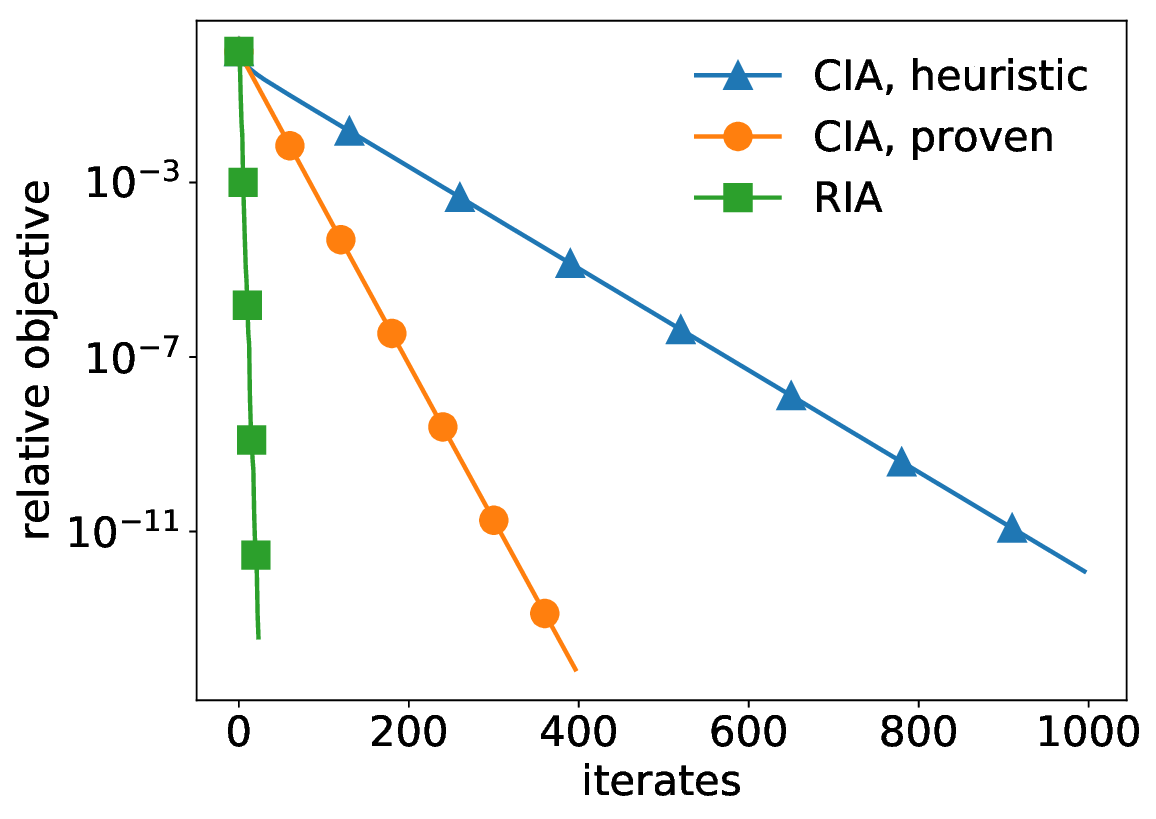}
	\caption{CIA and RIA methods applied to linear system, with matrix entries created from uniform distribution. CIA with the time step \(\tau = 1/[\sqrt{n} L]\) (\textbf{orange, circle}) performs better than the same method with heuristic time step \(\tau = 2/L\) (\textbf{blue, triangle}), but worse than RIA\@. This is the reverse of what was observed in previous examples.\label{fig:uniform}}
\end{figure}

The performance of coordinate descent methods and their optimal time steps vary significantly with the structure of the optimisation problem \citep{gur17, sun21, wri20}. If the linear systems above were constructed using a random distribution whose mean is not zero, then the results would look different. We demonstrate this with a numerical test with results in \figref{fig:uniform}.

We compare two time steps for the cyclic Itoh--Abe method, \(\tau_i = 2/L_i\) and \(\tau_i = 2/(L_i \sqrt{n})\), denoted by the curves labelled `heuristic' and `proven' respectively. While the heuristic time step was superior for most of the test problems considered in this section, it performs significantly worse for this example. Furthermore, in this case the randomised Itoh--Abe method converges faster than the cyclic one.

\subsection{Regularised logistic regression}

\begin{figure}[ht]
	\includegraphics[height=4cm]{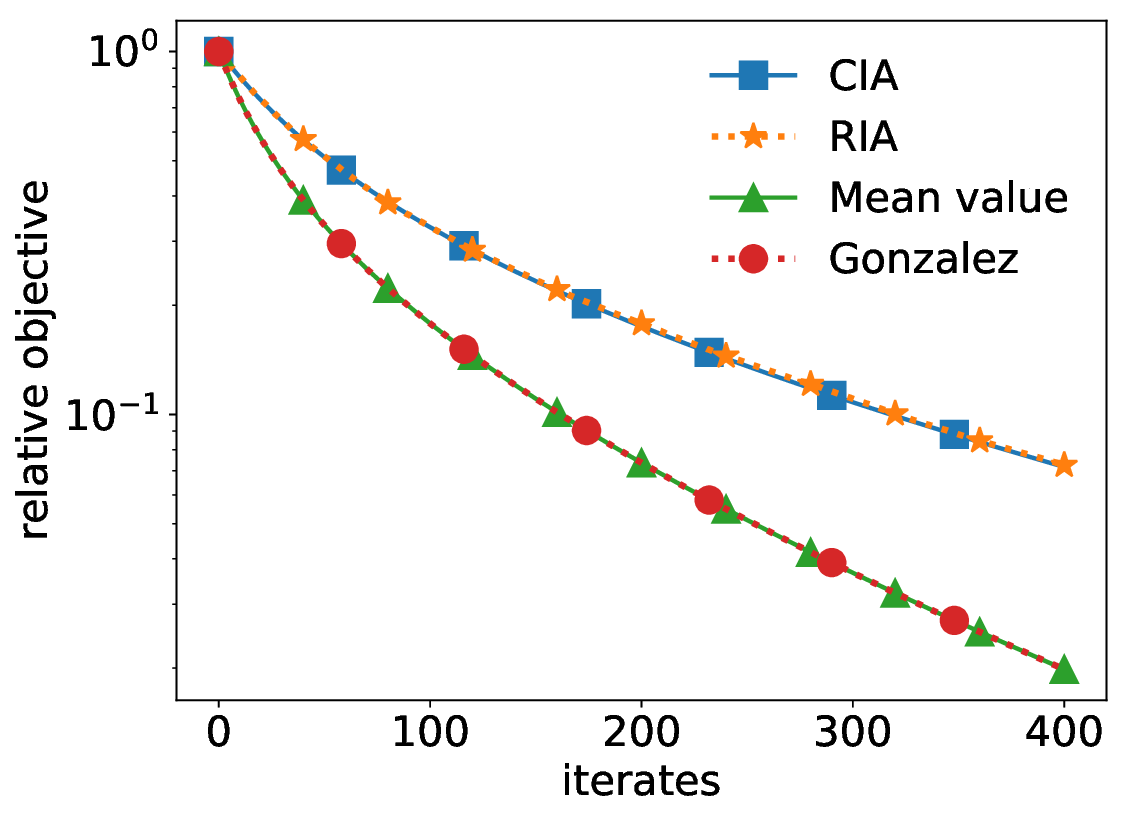}
	\caption{DG methods for \(l_2\)-regularised logistic regression~\eqref{eq:logreg_problem}. Convergence rate plotted as relative objective. The rates of randomised and cyclic Itoh--Abe methods almost coincide, and so do the mean value and Gonzalez discrete gradient methods.}\label{fig:logreg}
\end{figure}

We consider a \(l_2\)-regularised logistic regression problem, with training data \(\set{x^i, y_{i}}^{m}_{i=1}\), where \(x^i \in \RR^n\) is the data and \(y_i \in \{-1, 1\}\) is the class label. We wish to solve the optimisation problem
\begin{equation}\label{eq:logreg_problem}
	\min_{w \in \RR^n} V(w) = C \sum_{i=1}^m \log(1 + \mathrm{e}^{-y_i \inner{w, x^{i}}}) + \frac{1}{2} \|w\|^2,
\end{equation}
where \(C > 0\). We set \(n = 100\), \(m = 200\), \(C = 1\), and the elements of \({(y_{i})}_{i=1}^m\) is drawn from \(\{-1,1\}\) with equal probability. See \figref{fig:logreg} for the numerical results.

\subsection{Nonconvex function}

We solve the nonconvex problem
\begin{equation}\label{eq:nonconvex_problem}
	\min_{x \in \RR^n} V(x) = \|Ax\|^2 + 3\sin^2(\inner{c,x}),
\end{equation}
where \(A \in \RR^{n \times n}\) is a square, nonsingular matrix, and \(c\in \RR^n\) satisfies \(Ac = c\) and \(\|c\|=1\). This is a higher-dimensional extension of the scalar function \(x^2 + 3\sin^2(x)\) considered by~\cite{kar16}. This scalar function satisfies the P{\L} inequality~\eqref{eq:PL} for \(\mu = 1/32\), and it follows that \(V\) satisfies it for \(\mu = 1/(32\kappa)\), where \(\kappa\) is the condition number of \(A^*A\). Furthermore, the nonconvexity of \(V\) is observed by considering the restriction of \(V\) to \(x = \lambda c\) for \(\lambda \in \RR\), which has the form of the original scalar function. The function has the unique minimiser \(x^* = 0\).

\begin{figure}[ht]
	\centering
	\includegraphics[height=4cm]{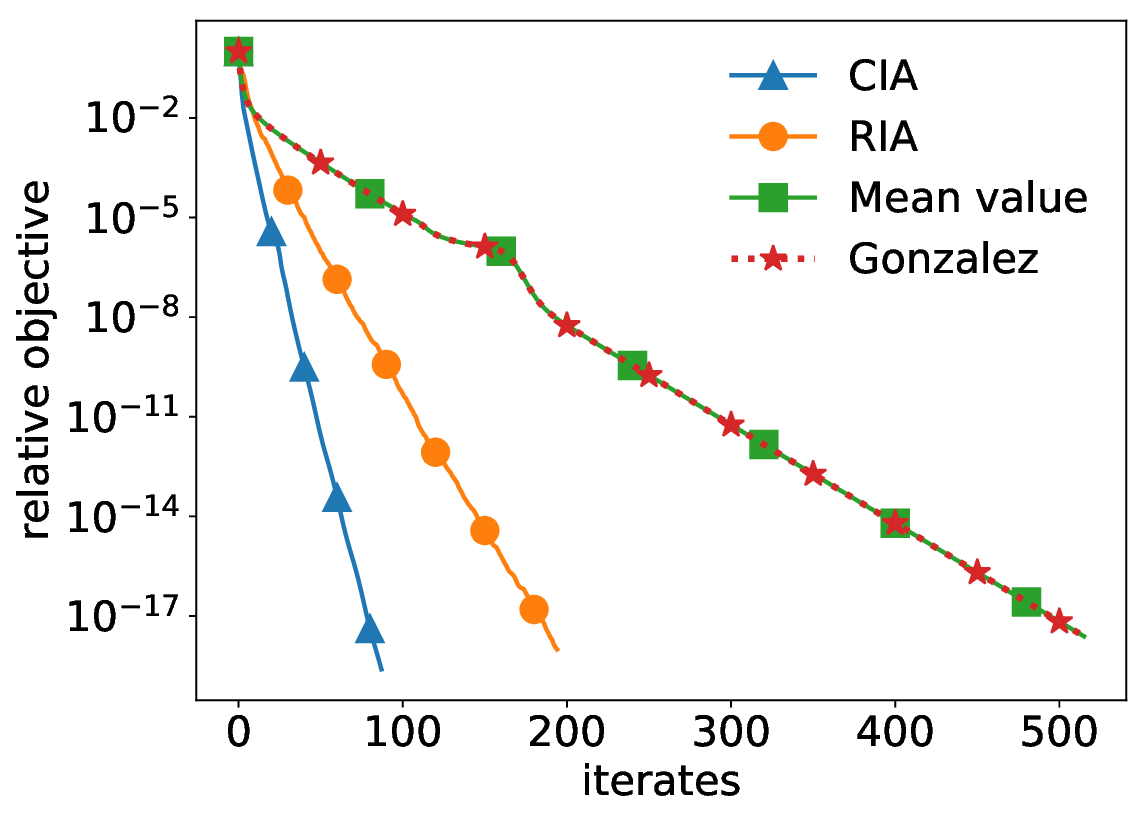}
	\includegraphics[height=4cm]{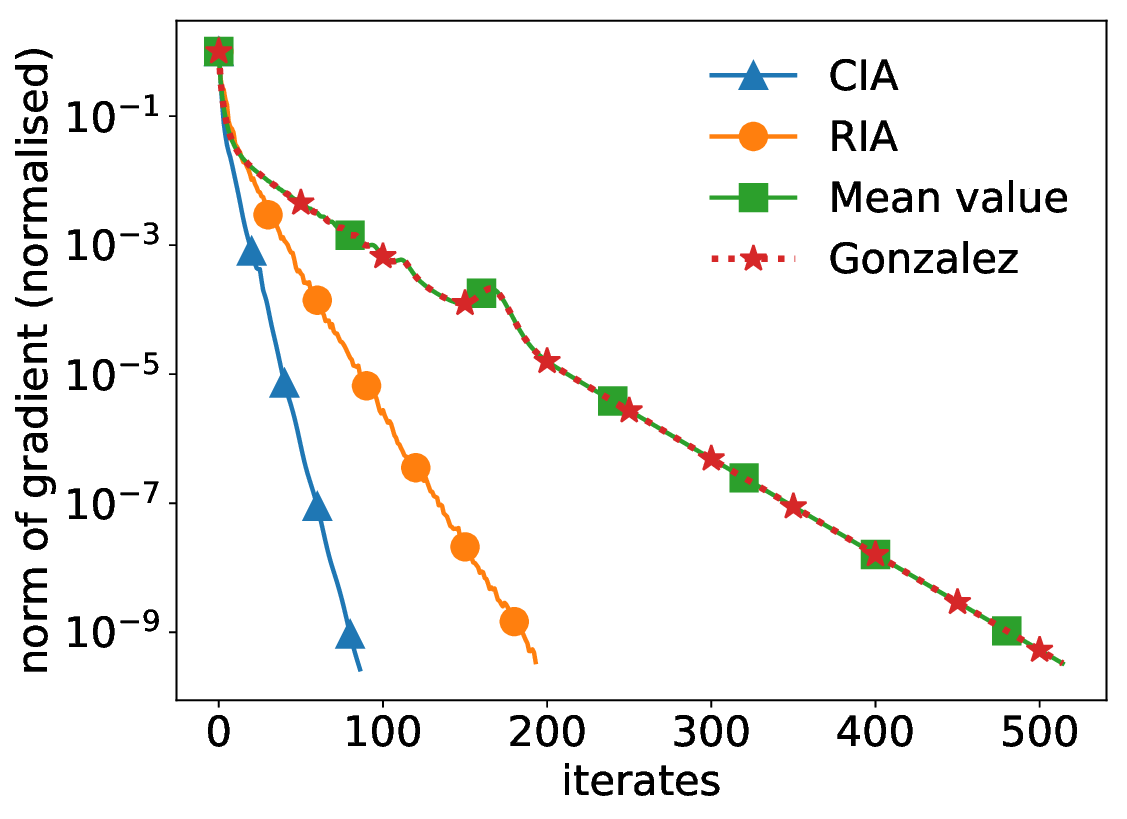}
	\caption{DG methods applied to nonconvex problem~\eqref{eq:nonconvex_problem} that satisfies the P{\L} inequality. \textbf{Left}: Relative objective. \textbf{Right}: Norm of gradient  \(\|\nabla V(x^{k})\| \; / \; \|\nabla V(x^0)\|\). Linear convergence rates are observed for the objective value and the gradient of the norm.}\label{fig:nonconvex}  % chktex 12
\end{figure}

We set \(n = 50\) and choose \(x^0\) constructed by random, independent draws from a Gaussian distribution. See \figref{fig:nonconvex} for the numerical results.

\subsection{Comparison of Itoh--Abe and coordinate descent for stiff problems}\label{sec:TV}

In image analysis and signal processing, variational optimisation problems often feature nonsmooth regularisation terms to promote sparsity, e.g.\ in the gradient domain or a wavelet basis. While one may replace these terms with smooth approximations, this can lead to stiffness of the optimisation problem, i.e.\ local, rapid variations in the gradient, requiring the use of severely small time steps for explicit numerical methods. In such cases, the cost of solving an implicit equation such as~\eqref{eq:dg_method} may be preferrable to explicit methods.

We investigate this scenario, by comparing the Itoh--Abe discrete gradient method to cyclic coordinate descent (CD)~\eqref{eq:CD} for solving (smoothed) total variation denoising problems. We denote by \(x^{\text{true}} \in \RR^n\) a ground truth image\footnote{We consider discretised images in two dimensions but express them in vector form for simplicity.}, and by \(x^\delta = x^{\text{true}} + \delta\) a noisy image, where \(\delta\) is random Gaussian noise. The total variation regulariser is defined as \(\TV(x) := \sum_{i=1}^n \|{[\nabla x]}_i\|\), with \(\nabla: \RR^n \to \RR^{ 2\times n}\) a discretised spatial gradient as defined in \citep{cha04}. As the nonsmoothness is induced by the \(\ell^2\)-norm, we approximate the regulariser by \(TV_\eps(x) := \sum_{i=1}^n \|{[\nabla x]}_i\|_{\eps}\), where \(\|x\|_\eps := \sqrt{\|x\|^2 + \eps}\).
The optimisation problem is thus given by
\begin{equation}\label{eq:total_variation}
	\argmin_{x \in \RR^n} \frac{1}{2} \|x - x^\delta\|^2 + \lambda \TV_\eps(x).
\end{equation}

Unless otherwise specified, the time step for CD is \(\tau_{\text{CD}} = 1/(2 \lambda \sqrt{\eps}+1)\) and for the Itoh--Abe discrete gradient method (DG) \(\tau_{\text{DG}} = 1/10\).

\begin{figure}[ht]
	\centering
	\includegraphics[height=4cm, clip, trim=1pt 0pt 30pt 0pt]{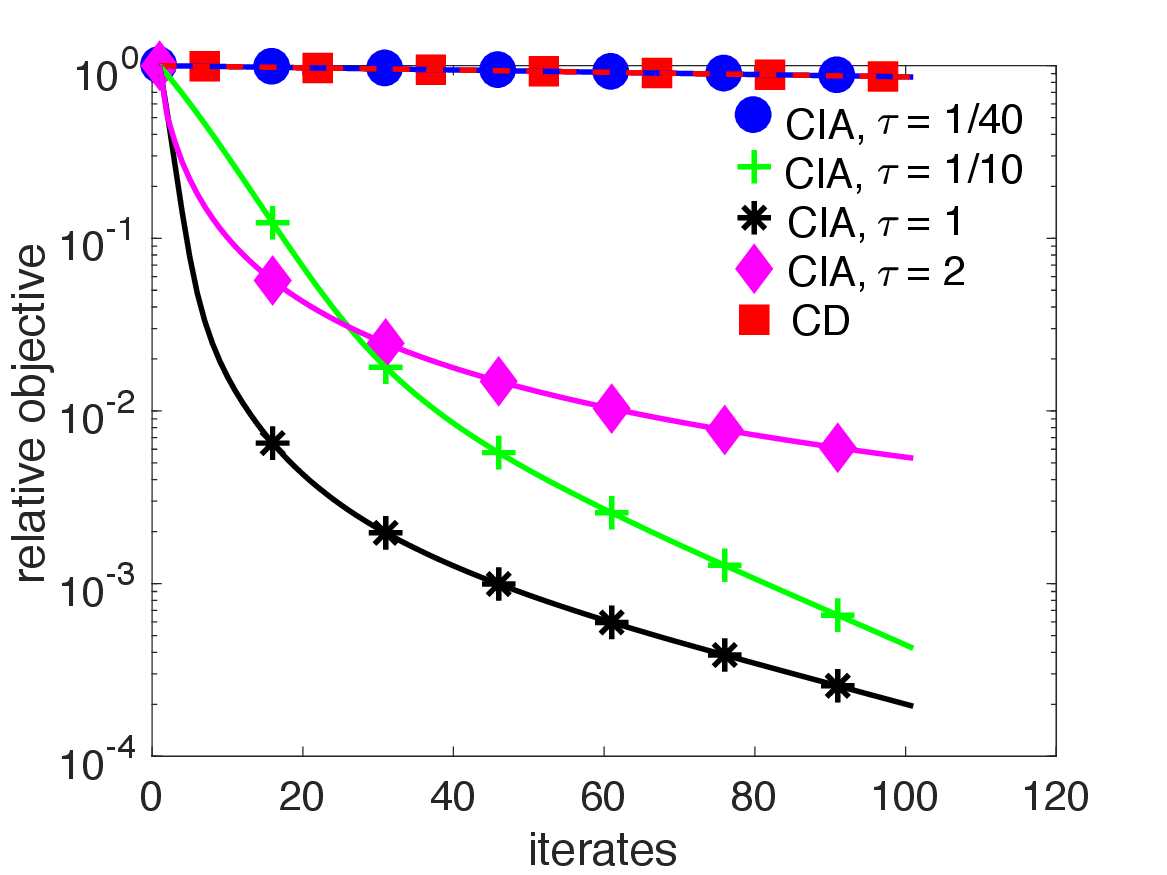}
	\includegraphics[height=4cm, clip, trim=100pt 0pt 80pt 0pt, width = 3.28cm]{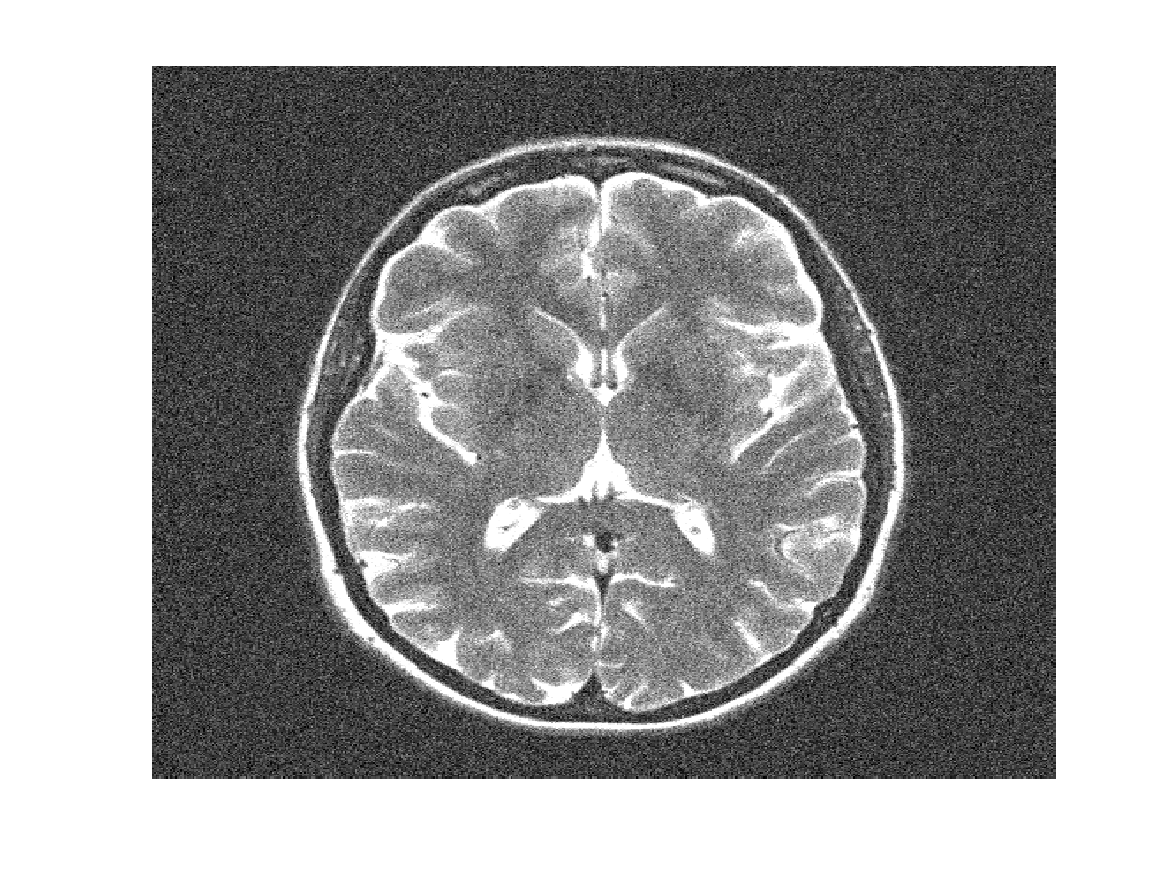}
	\includegraphics[height=4cm, clip, trim=100pt 0pt 80pt 0pt, width = 3.28cm]{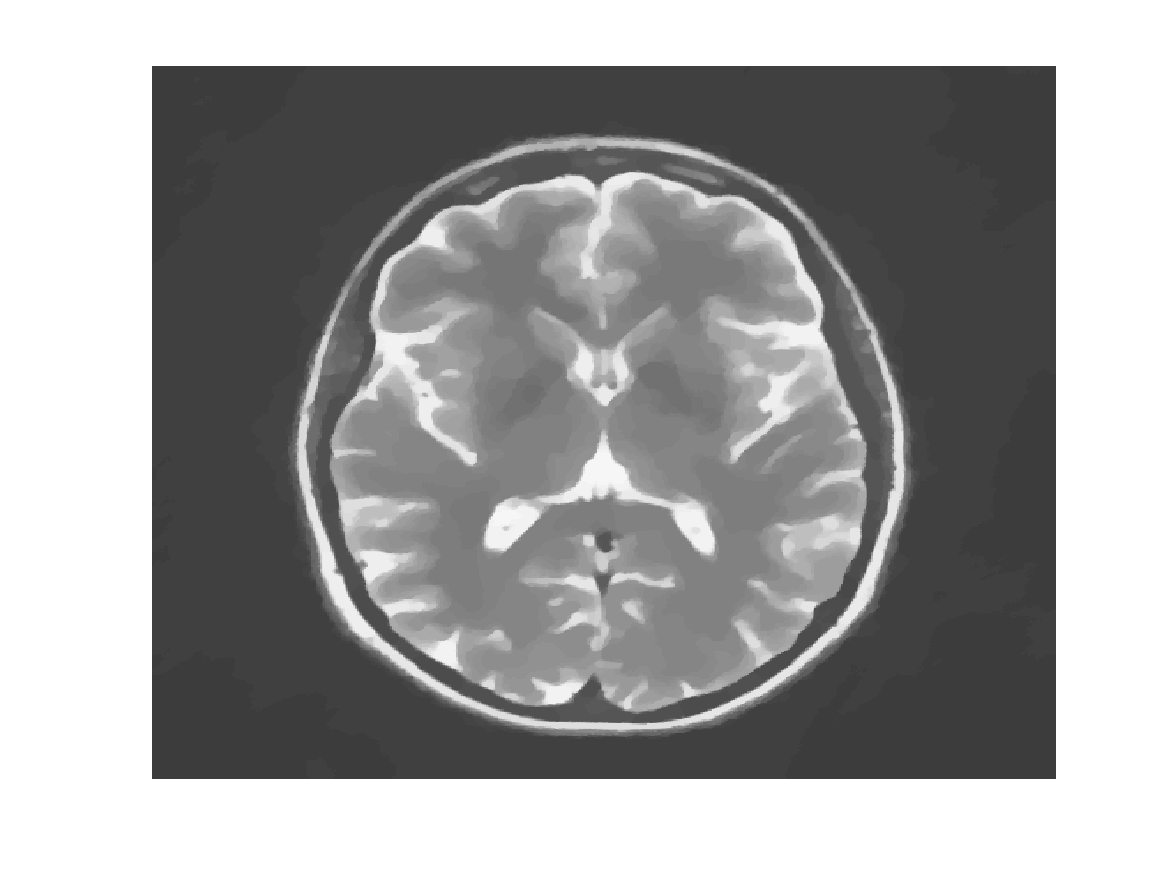}
	\caption{Comparison of CIA with different time steps to CD with standard time step \(1/L\) for solving the total variation problem~\eqref{eq:total_variation}. \textbf{Left}: Relative objective. Middle: Noisy image. Right: Reconstructed image. For the larger time steps, the CIA method converges significantly faster than both the CD method and the CIA method with the usual time step \(\tau_i = 2/L\).}\label{fig:dg_steps}  % chktex 12
\end{figure}

\begin{figure}[ht]
	\centering
	\includegraphics[height=3.6cm,clip, trim=7pt 0pt 40pt 0pt]{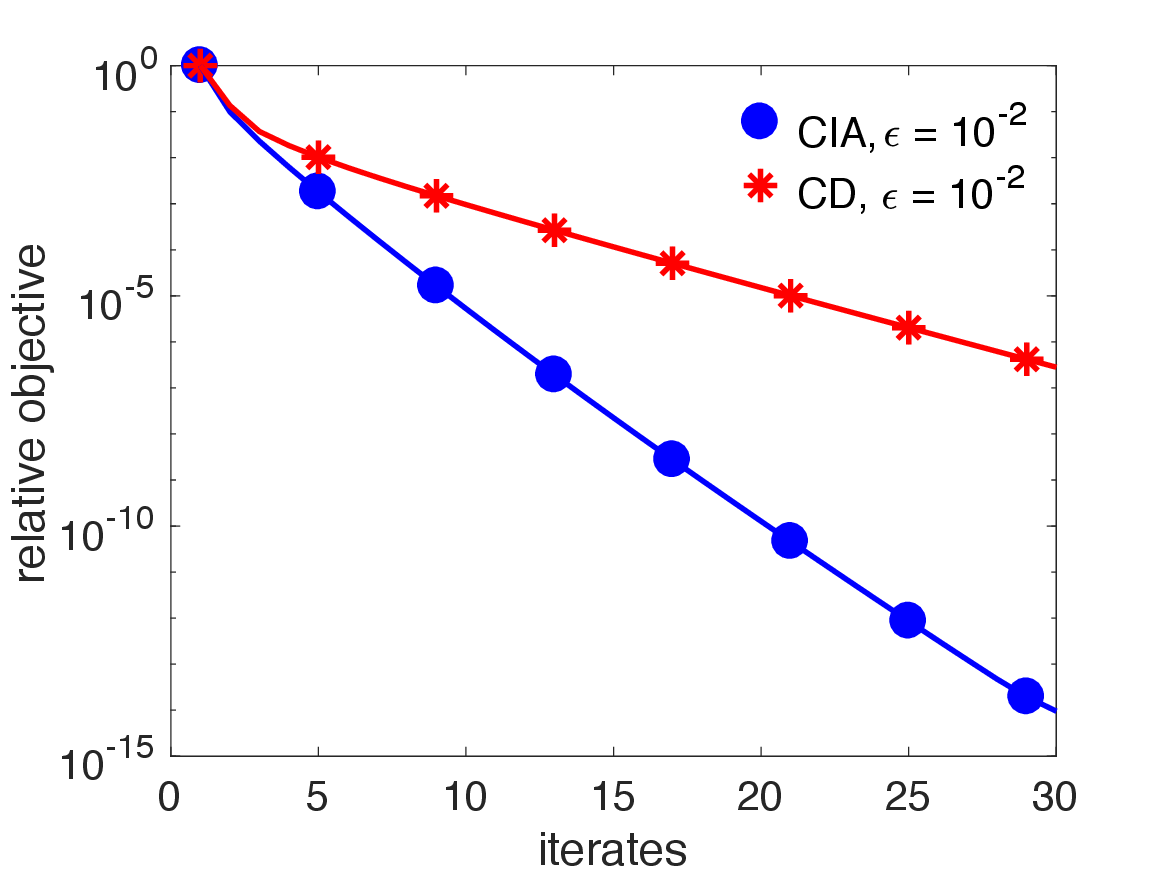}
	\includegraphics[height=3.6cm, clip, trim=30pt 0pt 40pt 0pt]{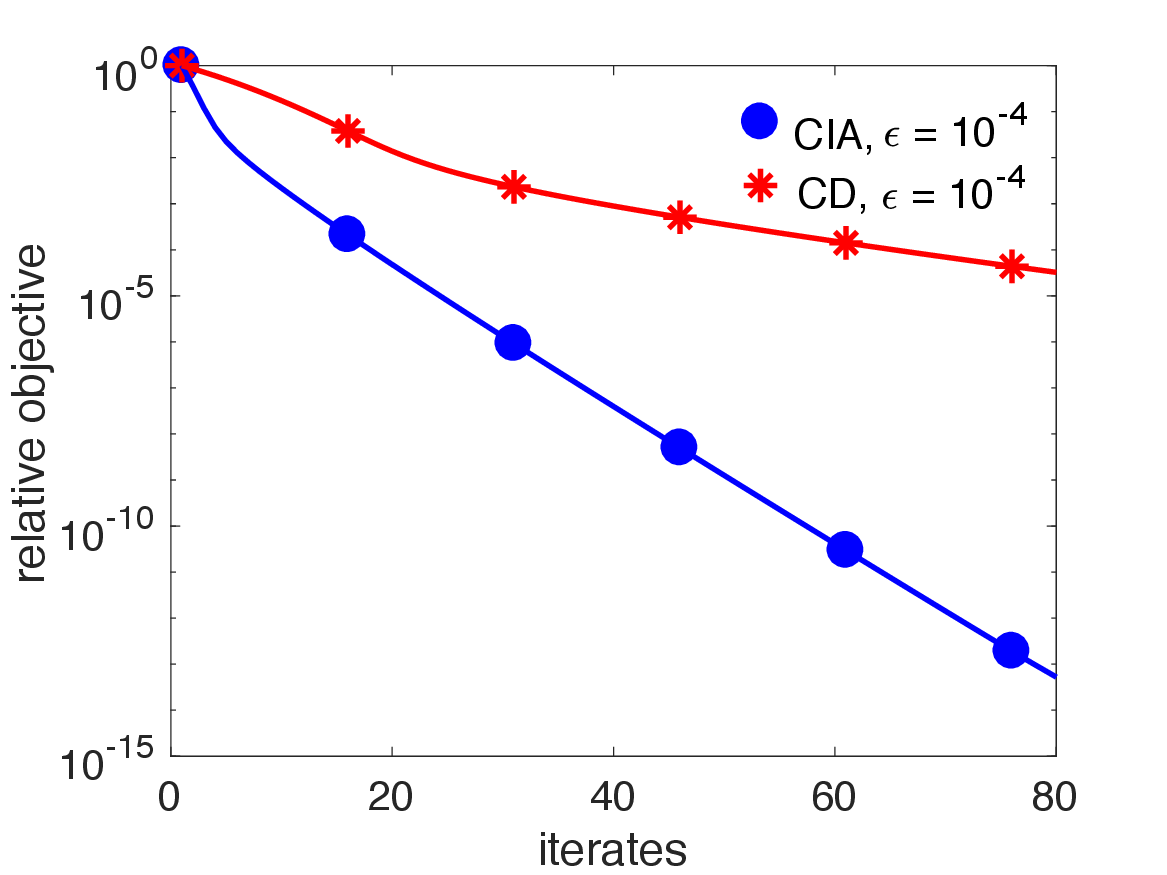}
	\includegraphics[height=3.6cm, clip, trim=30pt 0pt 40pt 0pt]{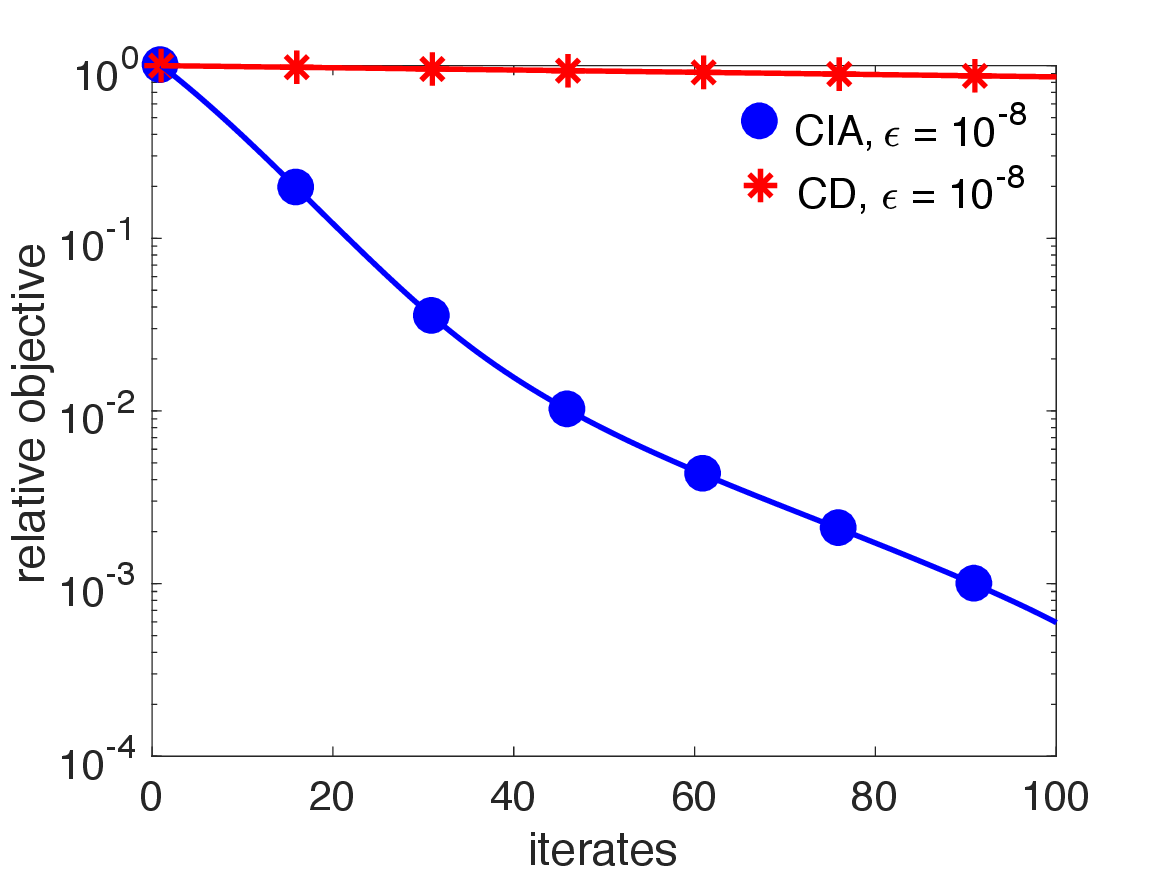}
	\caption{Comparison of CD to CIA for the total variation problem~\eqref{eq:total_variation} using three values of \(\eps\), (\textbf{top left}) \(10^{-2}\), (\textbf{top right}) \(10^{-4}\), and (\textbf{bottom}) \(10^{-8}\). The time steps are set to \(\tau_{\text{CIA}} = \sqrt{\tau_{\text{CD}}}\) where the latter time step is as usual.  As the problem becomes more ill-conditioned---when \(\eps\) is reduced---the performance of CIA improves relative to CD, demonstrating the CIA scheme's resilience to nonsmooth features.}\label{fig:epsilon}
\end{figure}

\begin{figure}[ht]
	\includegraphics[height=4cm]{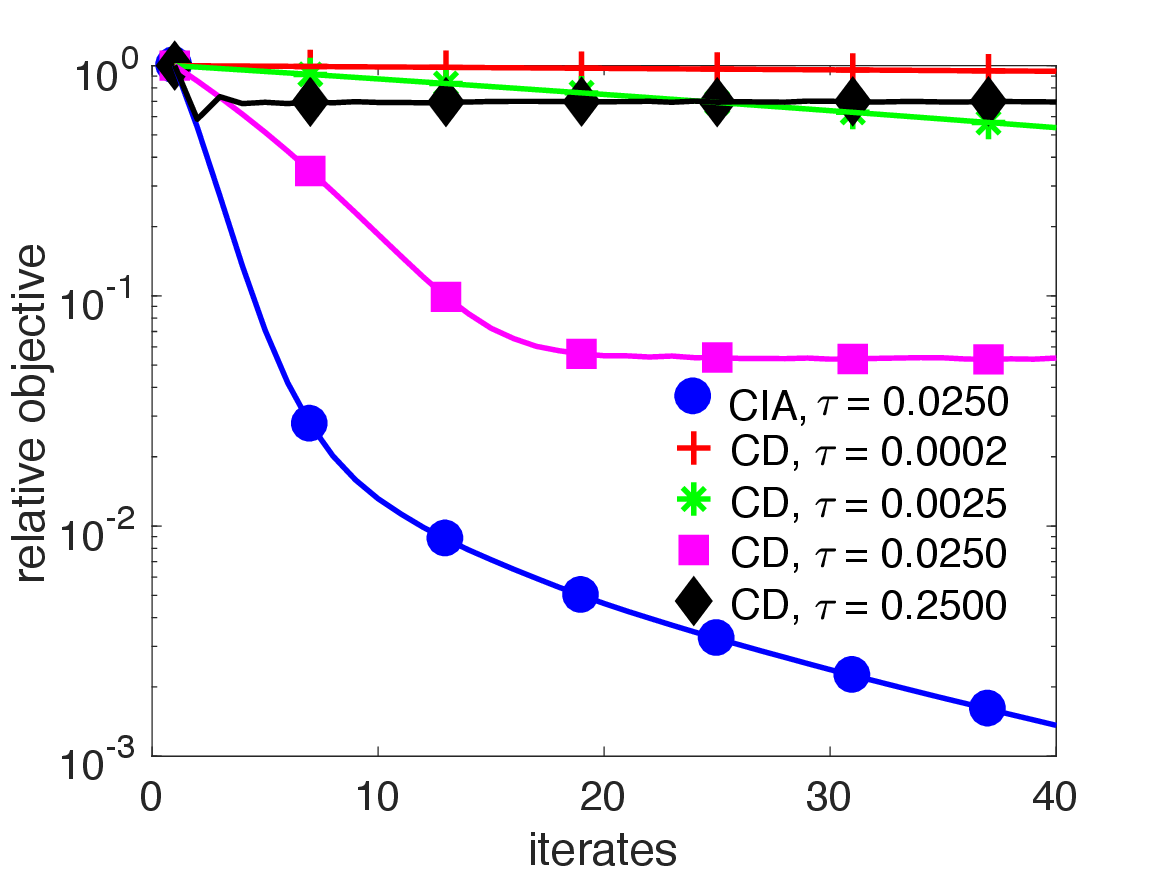}
	\caption{Comparison CD using different time steps versus CIA using a fixed time step, for the total variation problem~\eqref{eq:total_variation}. For smaller time steps, the CD iterates decrease too slowly, and for larger steps, they become unstable and fail to decrease.}\label{fig:cd_steps}
\end{figure}

\begin{figure}[ht]
	\includegraphics[height=4cm]{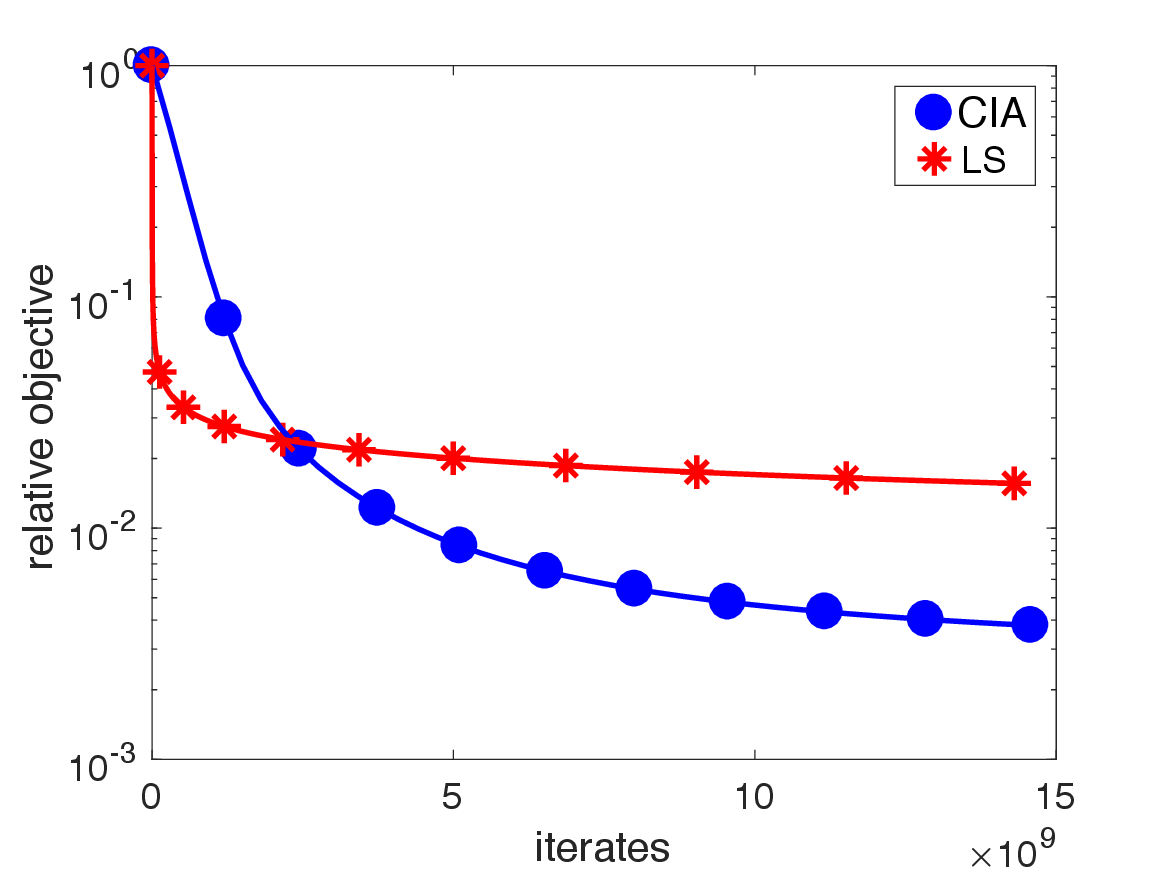}
	\caption{Comparison of CIA to backtracking line search for solving the total variation problem~\eqref{eq:total_variation} in terms of coordinate evaluations.}\label{fig:linesearch}
\end{figure}

In \figref{fig:dg_steps}, we compare the DG method for a range of time steps to CD\@. This demonstrates that the superior convergence rate of the DG method is stable with respect to a wide range of time steps. In \figref{fig:epsilon}, we compare the DG method to CD for different values of \(\eps\), demonstrating that the benefits of using the DG method increases as \(\eps\) gets smaller. In \figref{fig:cd_steps}, we compare different time steps for CD to the DG method, showing that for large time steps, the CD scheme is unstable and fails to decrease while for small time steps, the iterates decrease too slowly. In \figref{fig:linesearch}, we employ a simple backtracking line search (LS) method based on the Armijo-Goldstein condition, and compare this to the DG method.

\subsection{Comparison of methods for solving the discrete gradient equation}\label{sec:numerical_fixed_point}

We test the numerical performance of four methods for solving the discrete gradient equation~\eqref{eq:dg_method}, building on the fixed point theory in \secref{sec:fixed_point}.

The first method, denoted \textbf{F}, is the fixed point updates~\eqref{eq:fixed_point} proposed by~\cite{nor14} (\(\theta = 1\)). The second method, denoted \textbf{R}, is the relaxed fixed point method~\eqref{eq:relaxed}, where \(\theta\) is optimised according to~\eqref{eq:theta} if \(V\) is convex, and is otherwise set to \(1/2\). The third method, denoted \textbf{F+R}, is also the updates~\eqref{eq:relaxed} with \(\theta = 1\) by default, but whenever the discrepancy \(\|T(y^{k+1}) - y^{k+1}\|\) is greater than \(\|T(y^{k}) - y^k\|\), then the update is repeated with \(\theta\) set to half its previous value. This third option might be desirable in cases where \(\theta=1\) is expected to give faster convergence but also be unstable. The fourth method is the built-in solver \texttt{scipy.optimize.fsolve} in Python.

To test these methods, we performed 50 iterations of the discrete gradient method for different test problems, where at each iterate the discrete gradient solver would run until
\[
	\|r^k\|_\infty < \eps, \quad \mbox{where } r^k_i := \frac{y^k_i-y^{k-1}_{i}}{y^{k-1}_{i}} \mbox{ if } y^{k-1}_i \neq 0, \mbox{ and } r^k_i := y^k_i \mbox{ otherwise,}
\]
for some tolerance \(\eps > 0\), or until \(k\) reaches a maximum \(K_{\max}\). We then compare the average CPU time (\(s\)) for each of these methods. If a method fails to converge for a significant number of the iterations (\(> 10\%\)), we consider the method inapplicable for that test problem.

We test the methods for the mean value discrete gradient applied to three of the previous test problems, for \(\eps = 10^{-6}\) and \(10^{-12}\). We have not included results for the Gonzalez discrete gradient and other tolerances, as the results were largely the same.

The results are given in \tabref{tab:cpu}, where the best result is highlighted in bold for each test problem. We see that \textbf{R} is superior in stability, being the only method that locates the minimiser in every case. In all cases, \textbf{R} or \textbf{F+R} were the most efficient or close to the most efficient method. However, the relative performances of the different methods vary notably for the different test problems.

\begin{table}
	\caption{Average CPU time (\(s\)) over 50 iterations of~\eqref{eq:dg_method} with the mean value discrete gradient. Tolerance \(\eps = 10^{-6}\), \(10^{-12}\).\label{tab:cpu}}
	\begin{tabular}{@{}cccccc@{}}
		\toprule
		Test problem                                   & F              & R              & F + R          & fsolve & Tolerance    \\
		\midrule
		Linear system~\eqref{eq:linear_problem}        & N/A            & 0.006          & \textbf{0.002} & 0.190  & \(10^{-6}\)  \\
		Logistic regression~\eqref{eq:logreg_problem}  & \textbf{0.001} & 0.016          & 0.001          & N/A                   \\
		Nonconvex problem~\eqref{eq:nonconvex_problem} & N/A            & \textbf{0.003} & N/A            & N/A                   \\
		\bottomrule                                                                                                               \\
		\toprule
		Linear system~\eqref{eq:linear_problem}        & N/A            & 0.012          & \textbf{0.005} & 0.206  & \(10^{-12}\) \\
		Logistic regression~\eqref{eq:logreg_problem}  & 0.055          & 0.037          & \textbf{0.019} & N/A                   \\
		Nonconvex problem~\eqref{eq:nonconvex_problem} & N/A            & \textbf{0.005} & N/A            & 0.513                 \\
		\bottomrule % chktex 1
	\end{tabular}
\end{table}

\section{Conclusion}\label{sec:conclusion}

In this paper, we studied the discrete gradient method for optimisation, and provided several fundamental results on well-posedness, convergence rates and optimal time steps. We focused on four methods, using the Gonzalez discrete gradient, the mean value discrete gradient, the Itoh--Abe discrete gradient, and a randomised version of the Itoh--Abe method. Several of the proven convergence rates match the optimal rates of classical methods such as gradient descent and stochastic coordinate descent. For the Itoh--Abe discrete gradient method, the proven rates are better than previously established rates for comparable methods, i.e.\ cyclic coordinate descent methods \citep{wri15}.

There are open problems to be addressed in future work. Similar to acceleration for gradient descent, proximal gradient descent \citep{cha16}, and coordinate descent \citep{bec13, nes83, nes12, wri15}, we will study acceleration of the discrete gradient method to improve the convergence rate from \(\calO(1/k)\) to \(\calO(1/k^2)\). We would furthermore like to consider generalisations of the discrete gradient method to discretise gradient flows with respect to other measures of distance than the Euclidean inner product \citep{ben13, bur09}.

\section*{Acknowledgements}

All authors acknowledge support from CHiPS (Horizon 2020 RISE project grant). M. J. E., E. S. R., and C.-B. S. acknowledge support from the Cantab Capital Institute for the Mathematics of Information. M. J. E. and C.-B. S. acknowledge support from the Leverhulme Trust project on `
Breaking the non-convexity barrier', EPSRC grant Nr EP/M00483X/1, and the EPSRC Centre Nr EP/N014588/1.  Moreover, C.-B. S. acknowledges support from the RISE project NoMADS and the Alan Turing Institute.

\bibliographystyle{plainnat}
\bibliography{dg_refs_smooth_abbrev}

\appendix

%\section*{Appendix A.\ A measure-theoretical lemma}
%\label{app1}

\section{Bounds on discrete gradients}

\subsection{Proof of \lemref{lem:dg_assumption}}\label{sec:dg_bounds}

Part 1. We first consider the Gonzalez discrete gradient, and assume that \(x \neq y\) and write \(d := y-x\). Via the Gram-Schmidt process (see e.g.\ \citep[Section 4.6]{con94}), there is a vector \(d^{\perp}\) which satisfies \(\inner{d,d^\perp} = 0\), \(\|d^\perp\| = 1\), and
\[
	\overline{\nabla} V(x,y) =  \biginner{ \nabla V\del{\frac{x+y}{2}}, d^\perp} d^\perp + \frac{V(y) - V(x)}{\|y-x\|} d.
\]
By MVT, there is \(z \in [x,y]\) such that \(V(y) - V(x)  = \inner{\nabla V(z), y - x}\). Therefore, we obtain
\begin{equation}\label{eq:gonzalez_alt}
	\overline{\nabla} V(x,y) =  \biginner{ \nabla V\del{\frac{x+y}{2}}, d^\perp} d^\perp + \inner{\nabla V(z), d} d.
\end{equation}
From this, we derive \(\|\overline{\nabla} V(x,y)\|^2 \leq  \|\nabla V(\tfrac{x+y}{2})\|^2 + \|\nabla V(z)\|^2\). Thus property \emph{(i)} holds with \(C_n = \sqrt{2}\) and \(\delta \equiv 0\). To show property \emph{(ii)}, it is sufficient to note that since \(K\) is convex and has nonempty interior, then \(\nabla W\del{(x+y)/2} = \nabla V\del{(x+y)/2}\).

Part 2. Next we consider the mean value discrete gradient. It is clear that property \emph{(i)} holds with \(C_n = 1\) and \(\delta \equiv 0\). Property \emph{(ii)} is immediate from convexity of \(K\).

Part 3. For the Itoh--Abe discrete gradient, we set \(\delta(r) = r\). By applying the mean value theorem (MVT) to
\begin{equation}\label{eq:mvt}
	\del{\overline{\nabla} V(x, y)}_i = \frac{V(y_1, \ldots, y_i, x_{i+1}, \ldots, x_n) - V(y_1, \ldots, y_{i-1}, x_{i}, \ldots, x_n)}{y_i - x_{i}},
\end{equation}
we derive that \({(\overline{\nabla} V(x, y))}_i = \partial_i V(z^{i})\), where \(z^i = {[y_1, \ldots, y_{i-1}, c_i, x_{i+1}, \ldots, x_n]}^T\) for some \(c_i \in [x_i, y_i]\). Furthermore, we have \(\|z^i - x\| \leq \|y - x\|\), so \(z \in K_{\diam(K)}\). This implies that property \emph{(i)} holds with \(C_n = \sqrt{n}\). Property \emph{(ii)} is immediate.

\subsection{Proof of \lemref{lem:estimate_main}}\label{sec:estimate_main}

Part 1. Given \(z, d \in \RR^n\) as in~\eqref{eq:gonzalez_alt}, we denote by \(d^\perp\) a vector that satisfies
\[
	\inner{d, d^\perp} = 0, \qquad \|d^\perp\| =1, \qquad  \nabla V(x^{k}) = \inner{\nabla V(x^{k}) , d  }d + \inner{\nabla V(x^{k}), d^\perp} d^\perp.
\]
Note that by~\eqref{eq:dg_method}, it also holds that \(\inner{\overline{\nabla} V(x^k, x^{k+1}), d^\perp} = 0\). We compute
\begin{align*}
	\|\nabla V(x^{k})\|^2 & = \inner{\nabla V(x^{k}), d}^2 + \inner{\nabla V(x^{k}), d^\perp}^2                                                                       \\
	                      & \leq 2 \del{\|\overline{\nabla} V(x^k, x^{k+1})\|^2 +  \inner{\nabla V(x^{k}) - \nabla V(z), d}^2 + \frac{1}{4} L^2 \|x^k - x^{k+1}\|^2}.
\end{align*}
Since \(\inner{\nabla V(z), d} =(V(x^{k+1}) - V(x^{k}))/\|x^{k+1} - x^k\|\)
and \(d = \frac{x^{k+1} - x^{k}}{\|x^{k+1}-x^k\|}\), we have
\begin{align*}
	\inner{\nabla V(x^{k}) - \nabla V(z), d}^2 & = \frac{1}{\|x^k - x^{k+1}\|^2} \del{ \inner{\nabla V(x^{k}), x^{k+1} - x^{k}} - V(x^{k+1}) + V(x^{k}) }^2 \\ &  \leq \frac{1}{4} L^2 \|x^{k+1} - x^k\|^2,
\end{align*}
where the inequality follows from \propref{prop:smooth} \emph{(i)}. Applying~\eqref{eq:dissipation}, we conclude
\begin{align*}
	\|\nabla V(x^{k})\|^2 \leq 2 \del{\frac{1}{\tau_{k}} + \frac{1}{2} L^2 \tau_k } \del{ V(x^{k}) - V(x^{k+1})}.
\end{align*}

Part 2. We compute
\begin{align*}
	\|\nabla V(x^{k})\|^2 & \leq 2 \|\overline{\nabla} V(x^k, x^{k+1})\|^2 + 2 \lVert\int_0^1 \! \nabla V(sx^k + (1-s) x^{k+1}) - \nabla V(x^{k}) \dif s\rVert^2 \\
	                      & \leq 2 \|\overline{\nabla} V(x^k, x^{k+1})\|^2 + \frac{L^2}{2}  \|x^k - x^{k+1}\|^2                                                  \\& = 2 \del{ \frac{1}{\tau_{k}} + \frac{1}{4} L^2 \tau_k } \del{V(x^{k}) - V(x^{k+1})}.
\end{align*}

Part 3. We apply MVT like in~\eqref{eq:mvt} to obtain \(\del{\overline{\nabla} V(x^k, x^{k+1})}_i = \partial_i V(y^{i})\), where, for \(c_i \in [x^k_i, x^{k+1}_i]\),  \(y^i = {[x^{k+1}_1, \ldots, x^{k+1}_{i-1}, c_i, x^k_{i+1}, \ldots, x^k_n]}^T\). This gives
\begin{align*}
	\|\nabla V(x^{k})\|^2 & = \sum_{i=1}^n |\partial_i V(x^{k})|^2 \leq 2 \sum_{i=1}^n \del{ |\partial_i V(y^{i})|^2 + |\partial_i V(y^{i}) - \partial_i V(x^{k})|^2 } \\
	                      & \leq 2  \del{\|\overline{\nabla} V(x^k, x^{k+1})\|^2 + \overline{L}_{\Sum}^2 \|x^k - x^{k+1}\|^2}                                          \\& \leq 2\del{\frac{1}{\tau_{k}} + \overline{L}_{\Sum}^2 \tau_{k}} \del{V(x^{k}) - V(x^{k+1})}.
\end{align*}

Part 4. By respectively~\eqref{eq:lipschitz_direction} and~\eqref{eq:dissipation_IA}, we have
\begin{align*}
	\inner{\nabla V(x^{k}), x^k - x^{k+1}} & \leq V(x^{k}) - V(x^{k+1}) + \dfrac{L_{\max}}{2} \|x^k-x^{k+1}\|^2 \\& = \del{\dfrac{1}{\tau_{k}} + \dfrac{L_{\max}}{2}} \|x^k-x^{k+1}\|^2.
\end{align*}
Furthermore, \(\inner{\nabla V(x^{k}), x^k - x^{k+1}} = |\inner{\nabla V(x^{k}), d^{k+1}}| \|x^k - x^{k+1}\|\). From this, we derive
\begin{equation}\label{eq:random_1}
	\inner{\nabla V(x^{k}), d^{k+1}}^2 \leq \del{\frac{1}{\tau_{k}} + \frac{L_{\max}}{2}}^2 \|x^k - x^{k+1}\|^2.
\end{equation}
By the definition of \(\zeta\), we have
\begin{equation}\label{eq:random_2}
	\EE_{d^{k+1}}\inner{\nabla V(x^{k}), d^{k+1}}^2 \geq \zeta \|\nabla V(x^{k})\|^2.
\end{equation}
Combining~\eqref{eq:random_1} and~\eqref{eq:random_2}, we derive
\begin{align*}
	\|\nabla V(x^{k})\|^2 & \leq \frac{\tau_{k}}{\zeta} \del{\frac{1}{\tau_{k}} + \frac{L_{\max}}{2}}^2 \del{ V(x^{k}) - \EE_{d^{k+1}}[V(x^{k+1})]}.
\end{align*}
This concludes the proof.

\section{Convergence rate for cyclic coordinate descent}\label{sec:CD}

We now sharpen the convergence rates for cyclic coordinate descent (CD) \citep{bec13, wri15} to match those obtained for the Itoh--Abe discrete gradient method in \secref{sec:convergence_rate}. The CD method, for a starting point \(x^0\), time steps \(\tau_i > 0\), \(i = 1,\ldots, n\), and \(k = 0, 1, 2, \ldots\) is given by
\begin{equation}\label{eq:CD}
	\begin{aligned}
		x^{k,i+1} & = x^{k,i} - \tau_{i+1} \partial_{i+1} V(x^{k,i}) e^{i+1}, \qquad \mbox{for } i = 0, \ldots, n-1,
	\end{aligned}
\end{equation}
where \(x^{k,0} = x^k\) and \(x^{k+1} = x^{k,n}\). Recalling \secref{sec:convergence_rate}, we are interested in estimates for \(\beta > 0\) that satisfy~\eqref{eq:estimate_main}, where smaller \(\beta\) implies better convergence rate. In~\citep{bec13} (see Lemma 3.3) and referenced by~\cite{wri15}, the estimate
\begin{equation*}\label{eq:bec_optimal}
	\beta = 4 L_{\max}\del{1 + n L^2/L_{\min}^2},
\end{equation*}
is obtained, using the time step \(\tau_ i = 1/L_i\). This rate is optimised with respect to \(L_{\min}, L_{\max}\) when setting \(L_{\min} = L_{\max} = \sqrt{n} L\), yielding \(\beta = 8 \sqrt{n} L\). However, we show in \secref{sec:optimal_time_step} that the Itoh--Abe discrete gradient method achieves the stronger bound \(\beta = 4 \overline{L}_{\Sum} \leq 4 \sqrt{n} L\). We therefore include an analysis to show that the bound for CD can similarly be improved.

By the coordinate-wise descent lemma~\eqref{eq:lipschitz_direction}, we have
\begin{align*}
	V(x^{k,i})\! -\! V(x^{k,i+1}) & \geq  \inner{\nabla V(x^{k,i}), x^{k,i}-x^{k,i+1}} - \frac{L_{i}}{2}\|x^{k,i}-x^{k,i+1}\|^2 \\ & = (\tau_i - \frac{\tau_i^2 L_{i}}{2}) |\partial_{i+1} V(x^{k,i})|^2.
\end{align*}
For some \(\alpha \in (0,2)\), we choose \(\tau_i = \alpha/L_i\), and substitute into the above inequality to get
\begin{equation}\label{eq:descent1}
	V(x^{k,i}) - V(x^{k,i+1}) \geq \frac{\alpha}{L_{i}} \big(1 - \frac{\alpha}{2} \big) |\partial_{i+1} V(x^{k,i})|^2.
\end{equation}
We then compute
\begin{align*}
	\|\nabla V(x^{k})\|^2 & = \sum_{i=1}^n |\partial_i V(x^{k})|^2 \leq 2 \sum_{i=1}^n \del{|\partial_i V(x^{k}) - \partial_i V(x^{k,i-1})|^2 + |\partial_i V(x^{k,i-1})|^2 }                        \\
	                      & \overset{\eqref{eq:descent1}}{\leq} 2 \sum_{i=1}^n \del{  L^2 \|x^k - x^{k,i}\|^2 + \frac{L_{i}}{\alpha - \frac{\alpha^2}{2} } \del{V(x^{k,i-1}) - V(x^{k,i})} }         \\
	                      & \leq 2 \sum_{i=1}^n \del{  L^2 \sum_{j=0}^i \|x^{k,j} - x^{k,j+1}\|^2 +  \frac{L_{i}}{\alpha - \frac{\alpha^2}{2} }\del{V(x^{k,i-1}) - V(x^{k,i})} }                     \\
	                      & \leq 2\del{ \frac{n \alpha^2 L^2}{ L_{\min}^2} \sum_{j=0}^n |\partial_{j+1} V(x^{k,j})|^2 +  \frac{L_{\max}}{\alpha - \frac{\alpha^2}{2} } \del{V(x^{k}) - V(x^{k+1}) }} \\
	                      & \leq \frac{2 L_{\max} (1 + n \alpha^2 L^2/L_{\min}^2)}{\alpha - \frac{\alpha^2}{2} }  \del{V(x^{k}) - V(x^{k+1}) }.
\end{align*}
Setting \(\alpha = 1/\sqrt{n}\) and \(L_i = L\), we obtain the new estimate for \(\beta\),
\[
	\beta = 4 L \sqrt{n} \del{\frac{2 \sqrt{n}}{2 \sqrt{n} - 1}} \approx 4 \sqrt{n} L.
\]
This is approximately the same rate as that of the Itoh--Abe discrete gradient method.

Coordinate descent methods are often extended to \emph{block coordinate descent} methods. The above analysis can be extended to this setting by replacing \(n\) with the number of blocks.

\end{document}